\documentclass{article}
\usepackage[utf8]{inputenc}
\usepackage[margin=1in]{geometry} 
\usepackage{amsmath,amsthm,amssymb,amsfonts}
\usepackage{tikz}
\usetikzlibrary{shapes.geometric, arrows}
\usepackage{verbatim}
\usepackage{dsfont}
\usepackage{cite}
\usepackage{algpseudocode}
\usepackage{algorithm}
\usepackage{subcaption}
\usepackage[colorlinks=true]{hyperref}
\hypersetup{urlcolor=blue, citecolor=red, linkcolor=blue}
\usepackage[capitalise,noabbrev,nameinlink]{cleveref}
\usepackage[bottom]{footmisc}
\usepackage{multirow}
\usepackage[group-digits=false]{siunitx}
\usepackage{mwe}
\usepackage{mathtools}
\usepackage[inline]{enumitem}
\usepackage[]{mdframed}
\usepackage[skins,breakable]{tcolorbox}

\theoremstyle{plain}
\newtheorem{thm}{Theorem}[section] 
\newtheorem*{theorem*}{Main Result}
\newtheorem{lem}[thm]{Lemma}
\newtheorem{prop}[thm]{Proposition}

\newtheorem{cor}[thm]{Corollary}

\theoremstyle{definition}
\newtheorem{defn}[thm]{Definition}

\newtheorem{remark}[thm]{Remark}

\numberwithin{equation}{section}

\tikzstyle{dnc_step} = [rectangle, rounded corners, minimum width=3cm, minimum height=.7cm,text centered, text width = 4.5cm, draw=white]
\tikzstyle{arrow} = [draw, -latex',rounded corners]
\tikzstyle{point} = [circle, minimum width=.5cm, minimum height=.5cm, text centered, draw=black, fill=gray!30]
\tikzstyle{line} = [draw, -latex']

\usetikzlibrary{shapes,calc,arrows}

\begin{document}
\title{Structured Divide-and-Conquer for the Definite Generalized Eigenvalue Problem}
\author{James Demmel\footnote{Department of EECS (Computer Science Division) and Department of Mathematics, University of California Berkeley} \and Ioana Dumitriu\footnote{Department of Mathematics, University of California San Diego}
\and Ryan Schneider\footnote{Department of Mathematics, University of California Berkeley (ryan.schneider@berkeley.edu)}}
\date{}

\maketitle

\begin{abstract}
    This paper presents a fast, randomized divide-and-conquer algorithm for the definite generalized eigenvalue problem, which corresponds to pencils $(A,B)$ in which $A$ and $B$ are Hermitian and the Crawford number $\gamma(A,B) = \min_{\|x\|_2 = 1} |x^H(A+iB)x|$ is positive. Adapted from the fastest known method for diagonalizing arbitrary matrix pencils [Foundations of Computational Mathematics 2024], the algorithm is both inverse-free and highly parallel. As in the general case, randomization takes the form of perturbations applied to the input matrices, which regularize the problem for compatibility with fast, divide-and-conquer eigensolvers -- i.e., the now well-established phenomenon of \textit{pseudospectral shattering}. We demonstrate that this high-level approach to diagonalization can be executed in a structure-aware fashion by (1) extending pseudospectral shattering to definite pencils under structured perturbations (either random diagonal or sampled from the Gaussian unitary ensemble) and (2) formulating the divide-and-conquer procedure in a way that maintains definiteness. The result is a specialized solver whose complexity, when applied to definite pencils, is provably lower than that of general divide-and-conquer.
\end{abstract}

\small
\hspace{2mm} \textbf{Keywords:} Definite matrix pencil, Crawford number, generalized eigenvalue problem, Gaussian 

\hspace{20.7mm} unitary ensemble, pseudospectra \\

\vspace{-.3cm}
\indent \hspace{2mm} \textbf{MSC Class:} 15A22, 15B57, 65F15
\normalsize

\vspace{.3cm}



\section{Introduction}\label{section: intro}
\begin{figure}
    \centering
    \includegraphics[width=\linewidth]{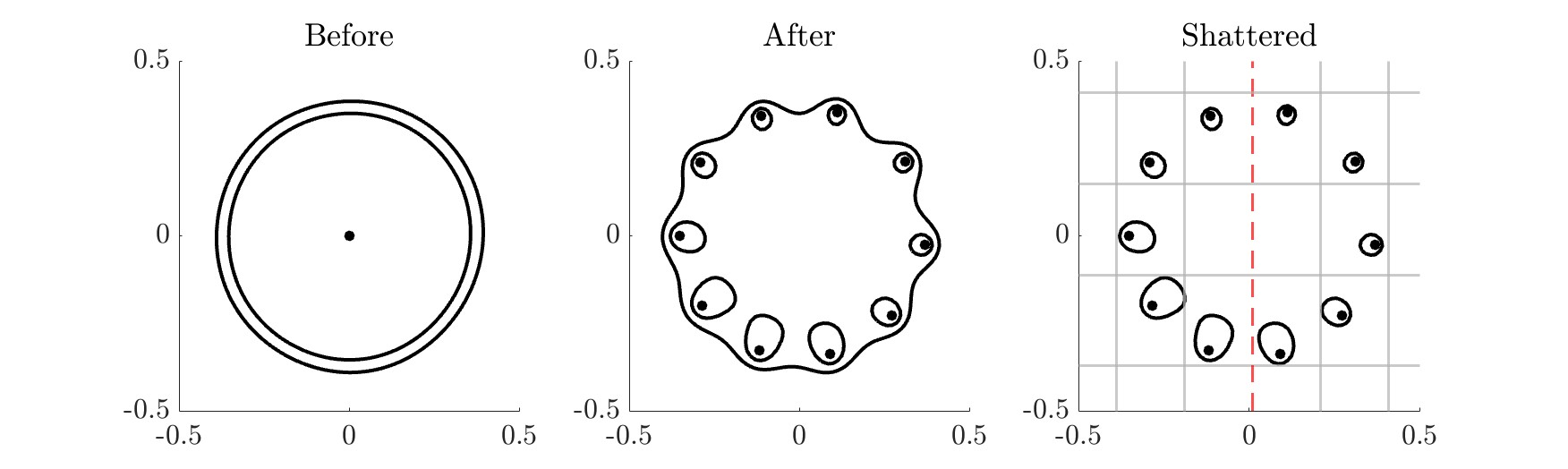}
    \caption{Pseudospectra of a $10 \times 10$ pencil $(A,B)$, constructed so that $B^{-1}A$ is a Jordan block,  before and after perturbation. The black curves trace the boundaries of $\epsilon$-pseudospectra for two choices of $\epsilon$; the eigenvalues they contained are marked with dark circles. In the rightmost plot, we include a random grid that shatters the tighter pseudospectrum of the perturbed problem, identifying in red a potential split for use in divide-and-conquer.}\label{fig: shattering_overview}
\end{figure}

The generalized eigenvalue problem seeks the eigenvalues and eigenvectors of the matrix pencil $(A,B) \in {\mathbb C}^{n \times n} \times {\mathbb C}^{n \times n}$. When $(A,B)$ is \textit{regular}, meaning the polynomial $p(x) = \det(A-xB)$ is not identically zero, its eigenpairs $(\lambda,v) \in {\mathbb C} \times {\mathbb C}^n\setminus \left\{ 0 \right\}$ can be defined as solutions to 
\begin{equation}\label{eqn: GEP}
    Av = \lambda Bv.
\end{equation}
If additionally $A$ and $B$ are Hermitian and the \textit{Crawford number}
\begin{equation}\label{eqn: crawfno}
    \gamma(A,B) \coloneqq \min_{\|x\|_2 = 1} |x^H(A+iB)x| = \min_{\|x\|_2 = 1} \sqrt{(x^HAx)^2 + (x^HBx)^2}
\end{equation}
is positive, then the pencil is \textit{definite}, in which case it is guaranteed to be regular with a full set of real eigenvalues. \\
\indent The definite generalized eigenvalue problem  appears frequently in applications, in part because it encompasses two important sub-problems:\ (1) the Hermitian, single-matrix eigenvalue problem, corresponding here to $B = I$, and (2) the generalized symmetric definite eigenvalue problem, in which $B$ and/or $A$ is positive definite. Often, the latter arises when $B$ (or $A$) is a Gram matrix -- as in, for example, extensions of classical support vector machines \cite{SVM} and quantum chemistry solvers \cite{FORD1974337}. Accordingly, specialized eigensolvers for this setting are well established. When $B$ is positive definite, solutions to \eqref{eqn: GEP} can be obtained by computing a Cholesky factorization $B = RR^H$ and diagonalizing the Hermitian matrix $R^{-1}AR^{-H}$ (e.g., \cite{Cholesky_symm_def}). Moreover, if we object to this approach on the basis of stability concerns,\footnote{Fast triangular inversion algorithms being only \textit{logarithmically} stable \cite{2007}.} we have plenty of alternatives to choose from \cite{Bunse84,Chandrasekaran00,Scott81,Veselié1993}. \\
\indent Of course, the pencil $(A,B)$ can be definite even if neither $A$ nor $B$ is positive definite itself. Indeed, noisy data or floating-point error may erase positive definiteness in the aforementioned examples without altering the definiteness of $(A,B)$. Unfortunately, algorithms that apply more broadly -- for example extensions of Jacobi's method that assume only symmetry \cite{Mehl_jacobi} -- cannot leverage the real spectra or stronger stability properties of definite pencils (see \cref{section: background}). We aim to fill this gap by tailoring fast, randomized, and highly parallel divide-and-conquer to arbitrary definite pencils, in particular without assuming that $A$ or $B$ is positive definite. \\
\indent As a high-level approach to eigenvalue problems, spectral divide-and-conquer has existed in the literature for decades \cite{BEAVERS1974143,bulgakov,MALYSHEV,Bai:CSD-94-793}. The specific algorithm we adapt, referred to as randomized or pseudospectral divide-and-conquer, was first introduced for the standard eigenvalue problem by Banks et al.\ \cite{banks2020pseudospectral} and subsequently extended to the generalized setting in \cite{arXiv}.  It is the fastest-known, fully general solver, capable of diagonalizing any matrix/pencil to accuracy $\varepsilon$ (in the spectral norm) in at most $O(n^{\omega_0} \log^2(n/\varepsilon))$ operations, where $O(n^{\omega_0})$ is the complexity of $n \times n$ matrix multiplication for $2 < \omega_0 \leq 3$. In the generalized case, divide-and-conquer aims to recursively divide $(A,B)$ into smaller pencils, whose eigenvalues form a disjoint union of the spectrum of $(A,B)$ and which can be diagonalized  independently and in parallel. Implementing this recursive division requires two key ingredients:\ (1) dividing lines/circles/curves that split eigenvalues into sets of roughly equal size and (2) efficient methods for computing spectral projectors onto the corresponding \textit{deflating subspaces} (see \cite{Projector_Paper} and the discussion in \cref{section: dnc}).  \\
\indent To obtain ingredient (1), the algorithms in \cite{banks2020pseudospectral} and \cite{arXiv} leverage the phenomenon of \textit{pseudospectral shattering}, wherein a small random perturbation to the input matrix/pencil guarantees, with high probability, both a minimum eigenvalue gap and minimally well-conditioned eigenvectors. This ultimately implies that the $\epsilon$-pseudospectrum of the perturbed problem is well-separated from a random grid containing its eigenvalues, provided $\epsilon$ is sufficiently small (see \cref{fig: shattering_overview}). In this context, the $\epsilon$-pseudospectrum of the pencil $(A,B)$ is defined as follows:
\begin{equation}\label{eqn: pseudospectrum}
    \Lambda_{\epsilon}(A,B) = \left\{z : \begin{array}{c}
        (A+E)u = z(B+F)u \\
        \text{for}  \; u \neq 0 \; \text{and} \; 
        \|E\|_2,\|F\|_2 \leq \epsilon  \\
        \end{array} \right\}.
\end{equation}
\indent Shattering suggests a straightforward procedure for diagonalizing $(A,B)$. First, add (small) random matrices to $A$ and $B$ to obtain a perturbed pencil $(\widetilde{A},\widetilde{B})$. Next, diagonalize $(\widetilde{A},\widetilde{B})$ by running divide-and-conquer over a random \textit{shattering grid}, whose boundary lines are well separated from a certain $\epsilon$-pseudospectrum of $(\widetilde{A},\widetilde{B})$ and can be efficiently searched over to find suitable eigenvalue splits. On the back end, the output of divide-and-conquer can be taken as an approximate diagonalization of $(A,B)$, which will be reasonably accurate (in the backward-error sense) if the initial perturbation is small. 

\indent Adapting this approach to the definite generalized eigenvalue problems requires addressing the following questions:
\begin{enumerate}
    \item Can we obtain a guarantee of pseudospectral shattering without destroying the structure of $(A,B)$? That is, can we ensure that the perturbed matrices $\widetilde{A}$ and $\widetilde{B}$ are both Hermitian with $(\widetilde{A},\widetilde{B})$ definite? Note that this cannot be done with the perturbations applied in \cite{arXiv}, which are Ginibre and therefore nonsymmetric in general. 
    \item If $(\widetilde{A},\widetilde{B})$ is definite -- and therefore has real eigenvalues -- can we accelerate some aspect of the divide-and-conquer procedure? In doing so, can we maintain definiteness through the recursive division?
\end{enumerate}

The answer to both questions is yes. To address the first, we prove a structured version of pseudospectral shattering under perturbations that are either (random) diagonal or sampled from the Gaussian unitary ensemble (see \cref{defn: GUE}). When such perturbations are sufficiently small, the perturbed pencil $(\widetilde{A},\widetilde{B})$ will remain definite, and in fact pseudospectral shattering can be defined entirely on the real line. This already implies efficiency gains in divide-and-conquer, swapping a two-dimensional shattering grid for a set of (real) splitting points.   \\
\indent Nevertheless, we can go further. By using a recently derived dynamically weighed Halley iteration to compute spectral projectors \cite{Projector_Paper}, we can reduce the asymptotic complexity of the algorithm, while maintaining definiteness, from that of \cite{arXiv} to $ O(n^{\omega_0} \log(n/\varepsilon)\log(\log(n/\epsilon) + \log(n)\log(\gamma^{-1})))$, where $\gamma$ is a lower bound on $\gamma(A,B)$ available on input. Such a lower bound is necessary to maintain definiteness in the algorithm; as we will see, both the allowable initial perturbation size as well as the parameters of divide-and-conquer will depend on $\gamma$.\footnote{Computational tools for obtaining this lower bound are discussed in \cref{section: background}.} As long as this lower bound is at least polynomial in $\varepsilon$ and $n^{-1}$, our specialized algorithm is provably faster than its predecessor, essentially trading a complexity factor of $\log(n/\varepsilon)$ for $\log(\log(n/\varepsilon))$. The key insight here is that spectral projectors can be computed much more efficiently if the corresponding eigenvalues are real.  \\
\indent For clarity, \cref{tab: operation_count} summarizes the main differences between the structured algorithm developed in this paper and the general version it adapts. Note that the respective accuracy guarantees hold with probability at least $1-O(n^{-1})$ for matrices satisfying $\|A\|_2, \|B\|_2 \leq 1$.

\renewcommand{\arraystretch}{1.5}
\begin{table}[h]
    \centering
    \begin{tabular}{c|p{2cm}p{2.6cm}p{3.2cm}p{4cm}}
         Version & Input & Output &  Guarantee & Complexity \\
         \hline
         General \cite{arXiv} & $A,B \in {\mathbb C}^{n \times n}$ $\varepsilon > 0$ & Invertible $S,T$ \; \; \; Diagonal $D$ & $||A-SDT^{-1}||_2 \leq \varepsilon $ \; \;  $||B - ST^{-1}||_2 \leq \varepsilon $ & $O(n^{\omega_0} \log^2(\frac{n}{\varepsilon}))$\\
         Structured  & $A,B \in {\mathbb C}^{n \times n}$ $\varepsilon > 0 \; \; \; \; \; \;  $ $\gamma \leq \gamma(A,B)$ & Invertible $X$ \; \; \; \; \; Diagonal $\Lambda_A,\Lambda_B$ & $||X^HAX - \Lambda_A||_2 \leq \varepsilon $ $||X^HBX - \Lambda_B||_2 \leq \varepsilon$ & $\begin{aligned}[t] O(&n^{\omega_0}\log(\textstyle\frac{n}{\varepsilon}) \log(\log(\textstyle\frac{n}{\varepsilon}) \;+ \\
         &\log(n)\log(\gamma^{-1})))
         \end{aligned} $\\
    \end{tabular}
\label{tab: operation_count}
\end{table}
\indent In an effort to balance stability and complexity, the general divide-and-conquer algorithm of \cite{arXiv} avoids matrix inversion. In large part, the decision to operate inverse free -- specifically when computing spectral projectors -- is rooted in the same stability concerns that advise against converting $(A,B)$ to the Hermitian problem $R^{-1}AR^{-H}$ (or more generally $B^{-1}A$), particularly when fast inversion algorithms are used.\footnote{See \cite[Chapter 6]{My_thesis} for a deeper discussion of precision requirements for divide-and-conquer with and without inversion.} At the same time, it also leaves open the possibility of minimizing associated (sequential or parallel) communication costs, where Ballard et al.\ established that inverse-free divide-and-conquer can attain communication lower bounds for general eigensolvers \cite{Ballard2010MinimizingCF}. Importantly, the modifications we apply to obtain our structured algorithm do not rely on inversion; hence, these advantages carry over.  \\
\indent The remainder of the paper is organized as follows. In \cref{section: background} we build motivation by presenting background information on definite pencils. We then prove our structured version of pseudospectral shattering in \cref{section: shattering} before developing the structured divide-and-conquer algorithm in \cref{section: dnc}. Finally, \cref{section: Examples} presents a collection of numerical examples to demonstrate the viability of structured divide-and-conquer in comparison to existing alternatives.

\subsection{Notation}
Throughout, $(A,B)$ is an $n \times n$ definite pencil with Crawford number $\gamma(A,B)$ and spectrum $\Lambda(A,B)$. $A^H$ and $A^{-H}$ denote the Hermitian transpose and inverse Hermitian transpose of $A$, respectively. $\|\cdot\|_2$ represents the spectral norm on matrices and the Euclidean norm on vectors, with $\kappa_2(\cdot)$ the spectral norm condition number and $\|(A,B)\|_2$ the spectral norm of the $n \times 2n$ block matrix $[A, \;B ]$. $\rho(A)$ is the spectral radius of $A$ while $\sigma_i(A)$ denotes its $i$-th singular value, with the standard convention $\|A\|_2 = \sigma_1(A) \geq \cdots \geq \sigma_n(A)$. When convenient, we may refer to the smallest singular value of $A$ as $\sigma_{\min}(A)$. Similarly, $\lambda_{\min}(A)$ denotes the smallest eigenvalue of the Hermitian matrix $A$. Finally, we use $\text{diag}(d_1,\ldots,d_n)$ to represent the (possibly block) diagonal matrix with diagonal entries (or blocks) $d_1, \ldots, d_n$.

\section{Background}\label{section: background}
In the numerical linear algebra literature, definite pencils were first considered rigorously by Crawford \cite{Crawford_1,Crawford_thesis}. In contrast to the generic setting, where a pencil may not even have a full set of eigenvalues, definite pencils are guaranteed to be both regular and diagonalizable. Moreover, their left/right eigenvectors are the same; hence, given a definite pencil $(A,B)$ there always exists an invertible matrix $X$ such that $X^HAX$ and $X^HBX$ are diagonal.
In the case where $A$ or $B$ is positive definite, the matrix $X$ can be obtained by way of a Cholesky factorization -- e.g., continuing the example from the introduction where $B = RR^H$, the unitary diagonalization $R^{-1}AR^{-H} = U\Lambda U^H$ implies that $(A,B)$ can be diagonalized by $X = R^{-H}U$. The eigenvalues in this case are clearly real, and this can similarly be shown for any definite pencil (see \cite[Theorem VI.1.18]{stewart1990matrix}). \\
\indent Note that the eigenvector matrix $X$ is \textit{not} unitary in general. Nevertheless, we can still expect stronger stability from the eigenvalues and eigenvectors of a definite pencil. This follows from a pair of results derived by Stewart, which are stated in terms of the chordal metric 
\begin{equation}\label{eqn: chordal_metrix}
    \chi(z,z') \coloneqq \frac{|z-z'|}{\sqrt{(|z|^2+1)(|z'|^2+1)}}.
\end{equation}

\begin{thm}[Stewart \cite{STEWART_DEFINITE}]\label{thm: definite_value_bound}
    Let $(A,B)$ be an $n \times n$ definite pencil and suppose that the Hermitian matrices $E,F \in {\mathbb C}^{n \times n}$ satisfy $\sqrt{\|E\|_2^2 + \|F\|_2^2} < \gamma(A,B)$.
    Then the pencil $(\widetilde{A}, \widetilde{B}) = (A+E,B+F)$ is definite. Moreover, if $\lambda_1 \leq  \cdots \leq \lambda_n$ and $\widetilde{\lambda}_1 \leq \cdots \leq \widetilde{\lambda}_n$ are the eigenvalues of $(A,B)$ and $(\widetilde{A},\widetilde{B})$, respectively, then for all $1 \leq i \leq n$,
    $$ \chi(\lambda_i,\widetilde{\lambda}_i) \leq \frac{\sqrt{\|E\|_2^2 + \|F\|_2^2}}{\gamma(A,B)}. $$ 
\end{thm} 

\begin{thm}[Stewart \cite{STEWART_DEFINITE}]\label{thm: definite_vector_bound}
    Let $(A,B)$ be an $n \times n$ definite pencil and suppose $(\widetilde{A}, \widetilde{B}) = (A+E,B+F)$ is also definite, where $E,F \in {\mathbb C}^{n \times n}$ are Hermitian. Let $v$ be a right eigenvector of $(A,B)$ corresponding to the eigenvalue $\lambda_1$. Suppose that $(\widetilde{A}, \widetilde{B})$ has eigenvalues $\widetilde{\lambda}_1, \ldots ,\widetilde{\lambda}_n$ and let $\delta = \min_{i > 1} \chi(\lambda_1, \widetilde{\lambda}_i)$. If $\sqrt{\|E\|_2^2 + \|F\|_2^2}< \delta \gamma(\widetilde{A}, \widetilde{B}), $
    then there exists a right eigenvector $\widetilde{v}$ of $(\widetilde{A}, \widetilde{B})$ corresponding to $\widetilde{\lambda}_1$ such that
    $$ \frac{\|\widetilde{v} - v\|_2}{\|v\|_2} \leq \frac{\sqrt{\|E\|_2^2 + \|F\|_2^2}}{\delta \gamma(\widetilde{A}, \widetilde{B})}. $$
\end{thm}

\indent \cref{thm: definite_value_bound} provides a criterion for ensuring that a perturbed pencil is definite. In fact, its proof implies 
\begin{equation}
    \aligned
    \gamma(\widetilde{A}, \widetilde{B})  
    & \geq \left[ 1 - \frac{\sqrt{\|E\|_2^2 + \|F\|_2^2}}{\gamma(A,B)} \right] \gamma(A,B).
    \endaligned 
    \label{eqn: Crawford_under_pertrub}
\end{equation}
This lower bound is tight:\ if $x_0$ is the unit vector at which the minimum in \eqref{eqn: crawfno} is achieved, it is easy to see that $(A-E,B-F)$ is indefinite for (rank one) perturbations $E = (x_0^HAx_0)x_0x_0^H$ and $F = (x_0^HBx_0)x_0x_0^H$, which by construction satisfy $\sqrt{\|E\|_2^2 + \|F\|_2^2} = \gamma(A,B)$. In this sense, $\gamma(A,B)$ can be interpreted as a distance to indefiniteness. \\ 
\indent Together, these results suggest the reciprocal of the Crawford number -- or more specifically the scale-invariant quantity $\|(A,B)\|_2/\gamma(A,B)$ -- as a condition number for the definite generalized eigenvalue problem. We can quantify this more precisely by considering the eigenvector matrix $X$ satisfying
\begin{equation}
    (X^HAX,X^HBX) = (\Lambda_A,\Lambda_B) = (\text{diag}(\alpha_1, \ldots, \alpha_n), \text{diag}(\beta_1,\ldots,\beta_n))
    \label{eqn: standard_choice_X}
\end{equation}
for $\alpha_i,\beta_i \in {\mathbb R}$ with $\alpha_i^2 + \beta_i^2 = 1$, which can be obtained from any diagonalizing matrix by scaling its columns accordingly. Importantly, the norms of $X$ and $X^{-1}$ are linked to $\gamma(A,B)$, as captured by the following observation of Elsner and Sun \cite[Proof of Theorem 2.3]{ELSNER1982341}. 
\begin{lem}[Elsner and Sun \cite{ELSNER1982341}]\label{lem: gamma_eigenvector}
Let $(A,B)$ be a definite pencil and let $X$ be a nonsingular eigenvector matrix satisfying \eqref{eqn: standard_choice_X}. Then 
$$\|X\|_2^2 \leq \gamma(A,B)^{-1} \; \; \text{and} \; \; \; \; \|X^{-1}\|_2^2 \leq \gamma(A,B)^{-1} \| (A,B)\|_2^2 $$
and therefore $ \kappa_2(X) \leq \frac{\|(A,B)\|_2}{\gamma(A,B)}. $
\end{lem}
\indent With this lemma in mind -- and specifically to obtain streamlined bounds -- we state subsequent results in terms of this quantity. It should be noted, however, that the bound implied by \cref{lem: gamma_eigenvector} is possibly quite loose. 

\subsection{Symmetric Pseudospectra}
We turn next to bounds on the pseudospectra of a definite pencil. Since we are interested in structured perturbations, we consider here a symmetrized version of \eqref{eqn: pseudospectrum}.\footnote{We are not the first to consider structured pseudospectra -- see e.g., \cite{TrefethenEmbree+2020,GRAILLAT200668}.}

\begin{defn}\label{defn: symm_pseudospecturum}
    The \textit{symmetric $\epsilon$-pseudospectrum} of $(A,B)$ is
    $$\Lambda_{\epsilon}^{\text{sym}}(A,B) = \left\{z : \begin{array}{l}
    (A+E)u = z(B+F)u \; \; \text{for} \; \; u \neq 0 \; \; \text{and} \\
    E,F \; \text{Hermitian with} \; \sqrt{\|E\|_2^2 + \|F\|_2^2} \leq \epsilon \\
    \end{array} \right\}.$$
\end{defn}

The decision to define $\Lambda_{\epsilon}^{\text{sym}}(A,B)$ in terms of $\sqrt{\|E\|_2^2 + \|F\|_2^2}$ is rooted in compatibility with \cref{thm: definite_value_bound}, which implies $\Lambda_{\epsilon}^{\text{sym}}(A,B) \subseteq {\mathbb R}$ provided $\epsilon < \gamma(A,B)$. In this setting, we can characterize $\Lambda_{\epsilon}^{\text{sym}}(A,B)$ entirely by the smallest singular value of $A-zB$ (equivalently, its distance to singularity). The proof of this result, which we omit for brevity, can be found in the full-length version of this paper \cite{structured_dnc_arxiv}.

\begin{lem}\label{lem: alt_sym_pseudo_equiv}
    Let $(A,B)$ be a definite pencil. If $\epsilon < \gamma(A,B)$ then
    $$ \Lambda_{\epsilon}^{\text{\normalfont sym}}(A,B) = \left\{ z : \sigma_{\min}(A-zB) \leq \epsilon \sqrt{1+|z|^2} \right\} \cap {\mathbb R}.$$
\end{lem}

\indent We can obtain similar bounds in terms of the general pseudospectrum \eqref{eqn: pseudospectrum}, in particular
\begin{equation}\label{eqn: pseudo_equiv}
    \Lambda_{\epsilon/\sqrt{2}}(A,B) \cap {\mathbb R} \subseteq \Lambda_{\epsilon}^{\text{ sym}}(A,B) \subseteq  \Lambda_{\epsilon}(A,B) \cap {\mathbb R},
\end{equation}
again assuming $\epsilon < \gamma(A,B)$. Of course, we always have $\Lambda_{\epsilon}^{\text{sym}}(A,B) \subseteq \Lambda_{\epsilon}(A,B)$. \\
\indent For finer control, we can derive a Bauer-Fike-type result -- i.e., upper and lower containments of $\Lambda_{\epsilon}^{\text{sym}}(A,B)$ in terms of intervals around the eigenvalues of $(A,B)$. For analogous bounds in the general setting, see \cite[Appendix A]{arXiv}.

\begin{thm}[Definite Bauer-Fike] \label{thm: symm_BF}
    Let $(A,B)$ be an $n \times n$ definite pencil with eigenvalues $\lambda_1, \ldots, \lambda_n \in {\mathbb R}$. For  $\epsilon < \min \left\{\sigma_n(B), \gamma(A,B) \right\}$ and $1 \leq i \leq n$ set 
    $$ r_{\epsilon} = \frac{\epsilon}{\gamma(A,B)}\left[1 + \left(\frac{\|A\|_2 + \epsilon}{\sigma_n(B) - \epsilon}\right)^2\right] \; \; \; \text{and} \; \; \;
    r_i = \begin{cases} 
      \frac{1}{\|B\|_2}, & \text{if} \; \; A = 0 \\
      \max \left\{ \frac{1}{\|B\|_2}, \frac{|\lambda_i|}{\|A\|_2} \right\}, & \text{otherwise.}
   \end{cases} $$ 
    Then
    $$ \bigcup_{i=1}^n (\lambda_i - \epsilon r_i, \lambda_i + \epsilon r_i ) \subseteq \Lambda_{\epsilon}^{\text{\normalfont sym}}(A,B) \subseteq \bigcup_{i=1}^n (\lambda_i - r_\epsilon, \lambda_i + r_{\epsilon}).$$
\end{thm}
\begin{proof}
     Recall that $z \in \Lambda_{\epsilon}^{\text{sym}}(A,B)$ implies $(A+E)u = z(B+F)u$ for some nonzero $u \in {\mathbb C}^n$ and Hermitian $E,F \in {\mathbb C}^{n \times n}$. Consequently, 
    \begin{equation}\label{eqn: bound_norm_of_z}
        |z| = \frac{\|(A+E)u\|_2}{\|(B+F)u\|_2 } \leq \frac{\|A\|_2 + \|E\|_2}{\sigma_n(B)-\|F\|_2} \leq \frac{\|A\|_2+\epsilon}{\sigma_n(B)-\epsilon}.
    \end{equation} 
    The upper containment now follows from \cref{thm: definite_value_bound}, allowing \eqref{eqn: bound_norm_of_z} to upper bound both $|\lambda_i|$ and $|\widetilde{\lambda}_i|$.  To obtain the lower inclusion we note that the eigenvalues of $(A+ \tau B, B)$ and $(A+\tau A, B)$ -- which are $\lambda_i + \tau$ and $\lambda_i (1+\tau)$, respectively -- belong to $\Lambda_{\epsilon}^{\text{sym}}(A,B)$, provided $\tau$ is real with $|\tau| \leq \epsilon / \|B\|_2$ or $|\tau| \leq \epsilon / \|A\|_2$.
\end{proof}

Both bounds in \cref{thm: symm_BF} are potentially quite loose. For the upper containment, this is a result of \eqref{eqn: bound_norm_of_z}, which is never attainable if $\sqrt{||E||_2^2 + ||F||_2^2} \leq \epsilon.$ The lower containment, on the other hand, is loose when $\epsilon$ approaches $\sigma_{\min}(B)$, in which case $\Lambda_{\epsilon}^{\text{sym}}(A,B)$ allows perturbations that  make $B$ nearly singular (and result in large eigenvalues). While it is likely possible to derive tighter bounds than this, \cref{thm: symm_BF} is sufficient for our purposes, as we demonstrate below.

\subsection{The Crawford Number, Numerically}
To close this section, we briefly discuss numerical methods for computing the Crawford number of $(A,B)$, recalling that our divide-and-conquer algorithm will require not only a guarantee that the input pencil is definite but also a lower bound $\gamma \leq \gamma(A,B)$. 
While identifying definiteness can be done cheaply -- see for example the arc algorithm introduced by Crawford and Moon \cite{Crawford_Moon} and recently revisited by Guo, Higham, and Tisseur \cite{arg_alg} -- obtaining $\gamma(A,B)$ to high accuracy is fairly demanding. Numerical tools for the latter are typically based around an alternative characterization of $\gamma(A,B)$:
\begin{equation}\label{eqn: alt_crawford_def}
    \gamma(A,B) = \max \left\{ \max_{\theta \in [0,2\pi)} \lambda_{\min}(\cos(\theta)A + \sin(\theta)B), \; \; 0  \right\}.
\end{equation}
Examples include a geometric algorithm of Uhlig \cite{Uhlig} and a subspace-acceleration-based method of Kressner, Lu, and Vandereycken \cite{subspace_crawfno}, both of which require repeatedly computing $\lambda_{\min}(\cos(\theta)A+\sin(\theta)B)$ for judicious choices of $\theta$. 
Since these algorithms are fairly involved, we also note that a relatively cheap lower bound for $\gamma(A,B)$ can be obtained from the smallest eigenvalue of $A$ (or $B$) if one of these matrices is positive (or negative) definite. 

\section{Symmetric Pseudospectral Shattering}\label{section: shattering}
We now prove a version of pseudospectral shattering for $\Lambda_{\epsilon}^{\text{sym}}(A,B)$. As discussed above, this specialized version can be stated entirely on the real line as long as the random perturbations to $A$ and $B$ are sufficiently small. In this context, shattering is defined relative to a set of points $g = \left\{g_i \right\} \subseteq {\mathbb R}$ with $g_i < g_{i+1}$. We say that $\Lambda_{\epsilon}^{\text{sym}}(A,B)$ is shattered with respect to $g$ if the following hold:
\begin{enumerate}
    \item Each eigenvalue of $(A,B)$ lies in a unique interval $(g_i,g_{i+1})$.
    \item $\Lambda_{\epsilon}^{\text{sym}}(A,B) \cap g = \emptyset$.
\end{enumerate}
The key insight here is that random perturbations guarantee not only a full set of distinct eigenvalues but also a minimum eigenvalue gap, which we denote as 
\begin{equation}\label{eqn: gap}
    \text{gap}(A,B) \coloneqq\min_{\substack{\lambda, \mu \in \Lambda(A,B) \\ \lambda \neq \mu}} |\lambda - \mu|.
\end{equation}

\subsection{Preliminaries}

\indent We start by considering perturbation strategies. In the general case, shattering is obtained for $\Lambda_{\epsilon}(A,B)$, with high probability, by adding Ginibre random matrices to $A$ and $B$. The natural analog for definite pencils is perturbation matrices sampled from the Gaussian unitary ensemble $\text{GUE}(n)$. 

\begin{defn}\label{defn: GUE}
    The $n \times n$ \textit{Gaussian unitary ensemble} $\text{GUE}(n)$ consists of matrices of the form 
    $$ Z = \frac{G+G^H}{\sqrt{2}}$$
    for $G$ an $n \times n$ Ginibre random matrix (with i.i.d.\ entries sampled from ${\mathcal N}_{\mathbb C}(0,\frac{1}{n})$).
\end{defn}

\indent Each $Z \in \text{GUE}(n)$ can be interpreted as a symmetrized Ginibre random matrix. Accordingly, we can easily write down probabilistic bounds for both the norm of $Z \in \text{GUE}(n)$ and the smallest singular value of a Hermitian matrix under GUE perturbation. The latter is obtained as a corollary of a result of Aizenman et al.\ \cite{Aizenman17}.

\begin{lem}\label{lem: GUE_norm_bound}
If $Z \in \text{\normalfont GUE}(n)$ then
$${\mathbb P} \left[ \|Z\|_2 \geq 4 + \sqrt{2} \right] \leq 2e^{-n}.$$
\end{lem}
\begin{proof}
Write $Z = \frac{G+G^H}{\sqrt{2}}$ for $G$ an $n \times n$ Ginibre random matrix. Since $\|Z\|_2 \leq \sqrt{2}\|G\|_2$, we have
    \begin{equation}
        {\mathbb P} \left[\|Z\|_2 \geq  4 + \sqrt{2} \right] \leq {\mathbb P} \left[ \|G\|_2 \geq 2\sqrt{2} + 1 \right] .
        \label{eqn: GUE_norm_bound_one}
    \end{equation}
    The result now follows from bounds on the Ginibre ensemble -- i.e., \cite[Lemma 2.2]{banks2020gaussian}.
\end{proof}

\begin{prop}\label{prop: smallest_GUE_bound}
    Let $M \in {\mathbb C}^{n \times n}$ be Hermitian and let $Z \in \text{\normalfont GUE}(n)$. For any $\mu > 0$, $t \leq \frac{\mu}{n}$, and an absolute constant $C < \infty$,
    $$ {\mathbb P} \left[ \sigma_n(M + \mu Z) \leq t \right] \leq C \frac{nt}{\mu}. $$
\end{prop}
\begin{proof}
    This is an immediate consequence of \cite[Theorem 1]{Aizenman17}.
\end{proof}

Of course, perturbations sampled from $\text{GUE}(n)$ are not the only option. We might instead consider using random diagonal matrices, noting that, in practice, $A$ and $B$ are often not only Hermitian but also banded. 

\begin{defn}\label{defn: diag_rand}
    Let $\rho$ be a probability density on ${\mathbb R}$ that is \textit{bounded} -- i.e., 
    $$ \|\rho\|_{\infty} \coloneqq\sup_{x \in {\mathbb R}} \rho(x) < \infty. $$ The corresponding \textit{diagonal random matrix} $D^{\rho}$ has nonzero entries sampled (independently) according to $\rho$.
\end{defn}

In the diagonal case, the analog of \cref{prop: smallest_GUE_bound} is the following.

\begin{prop}\label{prop: smallest_diagonal_bound}
    Let $M \in {\mathbb C}^{n \times n}$ be Hermitian and let $D^{\rho}$ be a diagonal random matrix corresponding to a bounded probability density $\rho$ on ${\mathbb R}$. For any $\mu > 0$ and  $t \leq \frac{\mu}{n}$,
    $${\mathbb P}[\sigma_n(M+\mu D^{\rho}) \leq t] \leq 2\|\rho\|_{\infty} \cdot \frac{nt}{\mu}. $$
\end{prop}
\begin{proof}
    If $\sigma_n(M+\mu D^{\rho}) \leq t$ then $\frac{1}{\mu}M + D^{\rho}$ has at least one eigenvalue in $[-\frac{t}{\mu}, \frac{t}{\mu}]$. Hence, letting $Y$ be the number of eigenvalues of $\frac{1}{\mu}M + D^{\rho}$ in this interval, Markov's inequality implies 
    \begin{equation}\label{eqn: }
        {\mathbb P}[\sigma_n(M+\mu D^{\rho}) \leq t] \leq {\mathbb E}[Y].
    \end{equation}
    The result now follows from work of Wegner \cite{Wegner} (see \cite[Equation 1.2]{Aizenman17}).
\end{proof}

Importantly, both options satisfy the following key result, which, as we will see, is the foundation of symmetric pseudospectral shattering.

\begin{thm}[Key Result]\label{thm: definite_key}
    Let $A \in {\mathbb C}^{n \times n}$ be a Hermitian matrix and let $I \subset {\mathbb R}$ be an interval with length $|I|$. Suppose the random matrix $V$ satisfies one of the following:
    \begin{enumerate}[itemsep=0em]
        \item $V \in \text{\normalfont GUE}(n)$.
        \item $V = D^{\rho}$ for a bounded probability density $\rho$ on ${\mathbb R}$.
    \end{enumerate}
    In either case, there exists a constant $C < \infty$ (uniform in $A$ and $n$) such that
    $$ {\mathbb P} \left[ A + V \; \text{has at least two eigenvalues in } I \right] \leq C |I|^2n^2. $$
\end{thm}
\begin{proof}
    The diagonal case was proved by Minami \cite{Minami} and  subsequently generalized to $\text{GUE}(n)$ by Aizenman et al.\ \cite{Aizenman17}. When $V = D^{\rho}$, $C$ is proportional to $\|\rho\|_{\infty}^2$.
\end{proof}

The assumption that $\rho$ is bounded in the diagonal case makes intuitive sense here; we cannot allow $\rho$ to be a Dirac function, for example, in which case the eigenvalues of $A+D^{\rho}$ are simply the eigenvalues of $A$ shifted by a constant. Appropriately, the factor of $\|\rho\|_{\infty}^2$ in \cref{thm: definite_key} weakens the result when $\rho$ is close to having such a ``spike." Of course, it also means that this bound -- as well as the diagonal version of pseudospectral shattering based on it -- is vacuous in a floating-point setting, where the support of $\rho$ is necessarily a set of measure zero (i.e., the set of floating-point numbers) and therefore $\|\rho\|_{\infty}  = \infty$. Nevertheless, we conjecture that shattering is still attainable in this setting (see \cref{fig: shattering_comparison}).\\
\indent \cref{thm: definite_key} pairs naturally with the following technical lemma, which is the final building block we need to establish shattering.

\begin{lem}\label{lem: singular_value_interval_bound}
Let $(A,B)$ be a definite pencil. If $(A,B)$ has at least $j$ eigenvalues in the interval $(z_0 - r, z_0 + r)$ for $z_0 \in {\mathbb R}$ and $r > 0$
then 
$$ \sigma_{n-j+1}(A - z_0B) \leq \frac{r \|(A, B)\|_2^2}{\gamma(A,B)}. $$
\end{lem}
\begin{proof}
    Let $X$ be a nonsingular eigenvector matrix of $(A,B)$ satisfying \eqref{eqn: standard_choice_X}. Standard singular value inequalities imply
    \begin{equation}
        \aligned 
        \sigma_{n-j+1}(A-z_0B) &= \sigma_{n-j+1}(X^{-H}(\Lambda_A - z_0 \Lambda_B)X^{-1}) \\
        &\leq \|X^{-H}\|_2 \sigma_{n-j+1}(\Lambda_A - z_0 \Lambda_B) \|X^{-1}\|_2 \\
        &= \|X^{-1}\|_2^2 \sigma_{n-j+1}(\Lambda_A - z_0 \Lambda_B).
        \endaligned 
        \label{eqn: apply_sv_inequal}
    \end{equation}
    Since $(A,B)$ has $j$ eigenvalues in $(z_0 - r, z_0 + r)$, at least $j$ diagonal entries $\alpha_i - z_0 \beta_i$ of $\Lambda_A - z_0 \Lambda_B$ satisfy
    \begin{equation}
        |\alpha_i - z_0 \beta_i| = |\beta_i| \left| \frac{\alpha_i}{\beta_i} - z_0 \right| \leq |\beta_i|r \leq r,
        \label{eqn: diag_entry_magnitude}
    \end{equation}
   meaning $\sigma_{n-j+1}(\Lambda_A - z_0 \Lambda_B) \leq r$. Applying this to \eqref{eqn: apply_sv_inequal} alongside \cref{lem: gamma_eigenvector} yields the final bound. 
\end{proof}

\subsection{Shattering Two Ways}

We are now ready to prove the main contribution of this section:\ symmetric pseudospectral shattering for definite pencils under GUE and random diagonal perturbations (respectively, \cref{thm: no_scaling_shattering} and \cref{thm: diagonal_shattering}). \\
\indent Since these results are fairly technical, we emphasize at the outset the important consequence they share. In each case, shattering is proved for a grid $g$ and a symmetric pseudospectrum $\Lambda_{\epsilon}^{\text{sym}}(A,B)$, where both $1/\epsilon$ and the number of grid points $|g|$ are polynomial in $n$ and inverse polynomial in $\mu$ (the ``size" of the perturbation). In the next section, we show that the complexity of our spectral divide-and-conquer algorithm will depend (at most) logarithmically on these quantities. Hence, the results of this section can be interpreted as follows:\ both GUE and random diagonal perturbations are sufficient to regularize the spectrum/pseudospectrum of $(A,B)$ to ensure efficiency in divide-and-conquer. 

\subsubsection{Version One:\ GUE}

We start with the GUE case.

\begin{thm}\label{thm: no_scaling_shattering}
    Let $(A,B)$ be an $n \times n$ definite pencil with $\|A\|_2,\|B\|_2 \leq 1$ and let $(\widetilde{A}, \widetilde{B}) = (A + \mu Z_1, B + \mu Z_2)$
    for independent $Z_1,Z_2 \in \text{\normalfont GUE}(n)$ and $0 < \mu < \frac{\gamma(A,B)}{12 \sqrt{2}}$. For $\omega = \frac{\mu^3\gamma(A,B)^2}{n^5}$ choose $z_0 \in  ( - \frac{3n^2}{2\mu} - \omega, -\frac{3n^2}{2\mu})$ uniformly at random and construct the grid of points 
    $$ g = \left\{ z_0 + j \omega : 0 \leq j \leq \left\lceil \frac{3n^2}{\mu \omega} \right\rceil +1 \right\}. $$ 
   Then $\Lambda_{\epsilon}^{\text{\normalfont sym}}(\widetilde{A}, \widetilde{B})$ is shattered with respect to $g$ for $\epsilon = \frac{\mu^6 \gamma(A,B)^2}{6n^{11}}$ with probability at least $1-O(n^{-1})$.
\end{thm}
\begin{proof}
     We condition on the events $\|Z_1\|_2, \|Z_2\|_2 < 6$ and  $\sigma_n(\widetilde{B}) > \frac{\mu}{n^2}$, which imply the following:
     \begin{enumerate}
         \item $(\widetilde{A},\widetilde{B})$ is definite with $\gamma(\widetilde{A},\widetilde{B}) \geq \frac{1}{2}\gamma(A,B)$. 
         \item$\Lambda(\widetilde{A},\widetilde{B}) \subset (-\frac{3 n^2}{2\mu}, \frac{3n^2}{2\mu})$.
         \item ${\mathbb P}\left[ \text{gap}(\widetilde{A},\widetilde{B}) \leq \delta \right] = O(\frac{n^4\delta}{\mu^5}) $. 
     \end{enumerate}
     Item (1) is a consequence of \cref{thm: definite_value_bound} and \eqref{eqn: Crawford_under_pertrub} while (2) follows from
     \begin{equation}\label{eqn: bound_spectral_norm}
        \rho(\widetilde{B}^{-1}\widetilde{A}) \leq \|\widetilde{B}^{-1}\widetilde{A}\|_2 \leq \frac{\|\widetilde{A}\|_2}{\sigma_n(\widetilde{B})} \leq \frac{3n^2}{2\mu},
     \end{equation}
     noting that $\|\widetilde{A}\|_2 \leq \frac{3}{2}$ since $\|A\|_2 \leq 1 $ and $\mu \|Z_1\|_2 < \frac{1}{2}$.\footnote{This follows from the fact that $\gamma(A,B) \leq \sqrt{2}$ if $\|A\|_2,\|B\|_2 \leq 1$.} To show (3), let ${\mathcal N}$ be a set of equally spaced points in $(-\frac{3n^2}{2\mu}, \frac{3n^2}{2\mu})$, with spacing $\delta/2$, such that $(-\frac{3n^2}{2\mu}, \frac{3n^2}{2\mu}) \subseteq \bigcup_{y \in {\mathcal N}}(y-\frac{\delta}{2},y+\frac{\delta}{2})$. It is easy to see that $|{\mathcal N}| \leq \lceil \frac{6n^2}{\mu\delta} \rceil$. Since $\Lambda(\widetilde{A},\widetilde{B})$ is contained in $(-\frac{3n^2}{2\mu},\frac{3n^2}{2\mu})$ we have
        \begin{equation}
            \aligned 
            {\mathbb P} \left[ \text{gap}(\widetilde{A}, \widetilde{B}) \leq \delta \right] &\leq {\mathbb P} \left[ |\Lambda(\widetilde{A}, \widetilde{B}) \cap (y-\delta,y+\delta) | \geq 2 \; \text{for some} \;  y \in {\mathcal N}  \right] \\
            & \leq | {\mathcal N} |  \max_{y \in {\mathcal N}} {\mathbb P} \left[ | \Lambda(\widetilde{A}, \widetilde{B}) \cap (y-\delta,y+\delta)| \geq 2 \right].
            \endaligned 
            \label{eqn: convert_prob_2}
        \end{equation}
        Now for any $y \in {\mathcal N}$, $|\Lambda(\widetilde{A},\widetilde{B}) \cap (y-\delta, y+\delta)| \geq 2$ implies
        \begin{equation}
            \sigma_{n-1}(\widetilde{A} - y\widetilde{B}) \leq  \frac{\delta \|(\widetilde{A}, \widetilde{B})\|_2^2}{ \gamma(\widetilde{A}, \widetilde{B})}
            \label{eqn: bootstrap_to_sv_bound}
        \end{equation}
         by \cref{lem: singular_value_interval_bound}. Applying $\gamma(\widetilde{A},\widetilde{B}) \geq \frac{1}{2}\gamma(A,B)$ and $\|(\widetilde{A}, \widetilde{B})\|_2 \leq \|\widetilde{A}\|_2 + \|\widetilde{B}\|_2 \leq 3$, this becomes
         \begin{equation}
            \frac{18\delta}{\gamma(A,B)} \geq \sigma_{n-1}(\widetilde{A} - y \widetilde{B}) = \mu \sigma_{n-1} \left( \frac{1}{\mu}[A - y \widetilde{B} ] + Z_1 \right).
            \label{eqn: rearrange_sv_bound}
        \end{equation}
         Thus, $|\Lambda(\widetilde{A}, \widetilde{B}) \cap (y-\delta, y+\delta) | \geq 2$ implies $\sigma_{n-1}(M + Z_1) \leq \frac{18\delta}{\mu\gamma(A,B)}$ for $M = \frac{1}{\mu}[A - y\widetilde{B}]$. But $M + Z_1$ is Hermitian and the singular values of a Hermitian matrix are the absolute values of its eigenvalues, so this is equivalent to $M + Z_1$ having at least two eigenvalues in the interval $[ - \frac{18 \delta}{\mu \gamma(A,B)}, \frac{18\delta}{\mu \gamma(A,B)} ]$, which occurs with probability at most $O(\frac{n^2\delta^2}{\mu^2 \gamma(A,B)^2})$ by \cref{thm: definite_key}. Plugging this into the union bound \eqref{eqn: convert_prob_2}, we obtain 
         \begin{equation}
            {\mathbb P} \left[ \text{gap}(\widetilde{A}, \widetilde{B}) \leq \delta \right] = O \left( \frac{n^4 \delta}{\mu^3 \gamma(A,B)^2} \right).
        \end{equation}
       \indent So far, letting $E$ be the event that $\|Z_1\|_2,\|Z_2\|_2 < 6$ and $\sigma_n(\widetilde{B}) > \frac{\mu}{n^2}$, we have shown 
    \begin{equation}
        {\mathbb P} \left[ \text{gap}(\widetilde{A}, \widetilde{B}) > \delta, \; \Lambda(\widetilde{A}, \widetilde{B}) \subseteq \left(-\frac{3n^2}{2\mu},\frac{3n^2}{2\mu}\right) \; \middle| \; E \right] \geq 1 - O \left(\frac{n^4 \delta}{\mu^3 \gamma(A,B)^2} \right).
        \label{eqn: conditional_prob}
    \end{equation}
    Choosing $\delta = \frac{\mu^3\gamma(A,B)^2}{n^5}$, this implies $\text{gap}(\widetilde{A}, \widetilde{B}) > \omega$ and $\Lambda(\widetilde{A},\widetilde{B}) \subseteq (-\frac{3n^2}{2\mu},\frac{3n^2}{2\mu})$ with probability at least $1 - O(n^{-1})$. In this case, by construction, each eigenvalue of $(\widetilde{A}, \widetilde{B})$ belongs to a unique interval $(g_i, g_{i+1})$ for $g_i \in g$. Moreover, the eigenvalues $\lambda \in \Lambda(\widetilde{A},\widetilde{B})$ are well separated from the grid points with high probability:\ applying the same geometric argument made in the proof of \cite[Theorem 3.12]{arXiv}, we have
    \begin{equation}
        {\mathbb P} \left[ \min_{\lambda,i} |\lambda - g_i| \leq \frac{\omega}{2n^2}\right] \leq \frac{1}{n}.
        \label{eqn: minimimum_g_dist}
    \end{equation}
    Thus, a simple union bound implies that each eigenvalue of $(\widetilde{A}, \widetilde{B})$ is contained in a unique grid interval, and is at least $\frac{\omega}{2n^2}$-away from the nearest grid point, with probability at least $1 - O(n^{-1})$. \\
    \indent When this occurs, shattering is guaranteed as long as $\Lambda_{\epsilon}^{\text{sym}}(\widetilde{A},\widetilde{B})$ is contained in a union of intervals of length $\frac{\omega}{n^2}$ centered at the eigenvalues of $(\widetilde{A}, \widetilde{B})$. Appealing to our version of Bauer-Fike for definite pencils (\cref{thm: symm_BF}), we know that if $\epsilon < \min \left\{ \gamma(\widetilde{A}, \widetilde{B}), \sigma_n(\widetilde{B}) \right\}$ then $\Lambda_{\epsilon}^{\text{sym}}(\widetilde{A}, \widetilde{B})$ is contained in a union of intervals of half-length
    \begin{equation} \label{eqn: apply_BF}
        r_\epsilon = \frac{\epsilon}{\gamma(\widetilde{A},\widetilde{B})}\left[1 + \left(\frac{\|\widetilde{A}\|_2 + \epsilon}{\sigma_n(\widetilde{B}) - \epsilon}\right)^2\right] \leq \frac{2\epsilon}{\gamma(A,B)}\left[ 1 + \left( \frac{\frac{3}{2}+\epsilon}{\frac{\mu}{n^2} - \epsilon} \right)^2\right].
    \end{equation}
    Further assuming $\epsilon < \frac{\mu}{2n^2}$ so that $\frac{3/2 + \epsilon}{\mu/n^2 - \epsilon} < \frac{4n^2}{\mu}$ and recalling $12 \sqrt{2}\mu < \gamma(A,B)$ we conclude
    \begin{equation}
        r_{\epsilon} \leq \frac{34 n^4 \epsilon}{12 \sqrt{2} \mu^3} \leq \frac{3n^4\epsilon}{\mu^3}.
        \label{eqn: simplified_re}
    \end{equation}
    Thus, we achieve shattering by taking $\frac{3 n^4 \epsilon}{\mu^3} \leq \frac{\omega}{2n^2} = \frac{\mu^3 \gamma(A,B)^2}{2n^7}$ or equivalently $\epsilon \leq \frac{\mu^6 \gamma(A,B)^2}{6n^{11}}$, which we note satisfies all of our assumptions on $\epsilon$. \\
    \indent In summary, we have shown that for $\epsilon = \frac{\mu^6 \gamma(A,B)^2}{6 n^{11}}$
    \begin{equation}\label{eqn: conditional_prob_of_shattering}
        {\mathbb P} \left[ \Lambda_{\epsilon}^{\text{sym}}(\widetilde{A},  \widetilde{B}) \; \text{is shattered with respect to} \; g \; | \; E \right] \geq 1 - O\left(n^{-1} \right).
    \end{equation}
    Since \cref{lem: GUE_norm_bound} and \cref{prop: smallest_GUE_bound} imply
    \begin{equation}\label{eqn: union_bound_for_conditioning}
        {\mathbb P} \left[ \|Z_1\|_2,\|Z_2\|_2 < 6, \; \sigma_n(\widetilde{B}) > \frac{\mu}{n^2} \right] \geq 1 -\frac{C}{n} - 4e^{-n}
    \end{equation}
    for an absolute constant $C < \infty$, we complete the proof by applying Bayes' theorem.
\end{proof}

\subsubsection{Version Two:\ Diagonal}
Swapping \cref{prop: smallest_GUE_bound} for \cref{prop: smallest_diagonal_bound}, we can easily bootstrap \cref{thm: no_scaling_shattering} into a proof of shattering under diagonal perturbations, in particular for the same choices of $g$ and $\epsilon$. 

\begin{thm}\label{thm: diagonal_shattering}
    Let $\rho$ be a bounded probability density on ${\mathbb R}$ such that 
    $${\mathbb P} [|Y| > M_{\rho} \; | \; Y \sim \rho] \leq \theta $$
    for $M_{\rho} > 0 $ and $\theta \in (0,1)$. Given an $n \times n$ definite pencil $(A,B)$ with $\|A\|_2,\|B\|_2 \leq 1$, let $(\widetilde{A},\widetilde{B}) = (A+\mu D_1^{\rho}, B + \mu D_2^{\rho}) $
    for independent, diagonal random matrices $D_1^{\rho}$, $D_2^{\rho}$ and $0 < \mu < \frac{\gamma(A,B)}{2\sqrt{2}M_{\rho}}$. Then $\Lambda_{\epsilon}^{\text{\normalfont sym}}(\widetilde{A},\widetilde{B})$ is shattered with respect to $g$, for both $g \subseteq {\mathbb R}$ and $\epsilon > 0$ as in \cref{thm: no_scaling_shattering}, with probability at least 
    $$\left(1 - \frac{C \|\rho\|_{\infty}^2}{n} \right) \cdot\left( 1 - 2n \theta - \frac{2\|\rho\|_{\infty}}{n} \right),$$
    where $C < \infty$ is an absolute constant.
\end{thm}

\begin{proof}
    We condition on the events $\|D_1^{\rho}\|_2, \|D_2^{\rho}\|_2 < M_{\rho}$ and $\sigma_n(\widetilde{B}) > \frac{\mu}{n^2}$. Since $D_1^{\rho}$ is diagonal -- and therefore has norm equal to the magnitude of its largest diagonal entry -- a simple union bound implies ${\mathbb P}[\|D_1^{\rho} \|_2 \geq M_{\rho}] \leq n\theta$. Similarly, we have 
    \begin{equation}\label{eqn: conditional_one}
        {\mathbb P}\left[\sigma_n(\widetilde{B}) \leq \frac{\mu}{n^2} \right] \leq \frac{2 \|\rho\|_{\infty}}{n}
    \end{equation}
    by \cref{prop: smallest_diagonal_bound}. Hence, $\|D_1^{\rho}\|_2,\|D_2^{\rho}\|_2 < M_{\rho}$ and $\sigma_n(\widetilde{B}) > \frac{\mu}{n}$ occur simultaneously with probability at least $1 - 2n\theta - \frac{2 \|\rho\|_{\infty}}{n}$. 
    Moreover, they imply the following:
    \begin{enumerate}
        \item $\gamma(\widetilde{A},\widetilde{B}) \geq \left( 1 - \frac{\sqrt{2}\mu M_{\rho}}{\gamma(A,B)} \right) \gamma(A,B) \geq \frac{1}{2}\gamma(A,B)$ (by \eqref{eqn: Crawford_under_pertrub} and our bound on $\mu$).
        \item $\|\widetilde{A}\|_2 \leq \|A\|_2 + \mu M_\rho \leq \frac{3}{2}$ (since $\mu M_{\rho} \leq \frac{\gamma(A,B)}{2\sqrt{2}}$ and $\gamma(A,B) \leq \sqrt{2}$). 
    \end{enumerate}
    The proof of shattering from \cref{thm: no_scaling_shattering} can now be repeated directly. In this case, the probability that $\text{gap}(\widetilde{A},\widetilde{B}) < \delta$ is accompanied by a constant proportional to $\|\rho\|_{\infty}^2$ (via \cref{thm: definite_key}). To highlight its dependence on $\rho$, we can rewrite the success probability as $\left(1 - \frac{C \|\rho\|_{\infty}^2}{n} \right) \cdot \left( 1 - 2n\theta - \frac{2 \|\rho\|_{\infty}}{n} \right)$ for an absolute constant $C < \infty$.
    \end{proof} 
    
 \indent Here, the probabilistic bound ${\mathbb P}[|Y| > M_\rho \; | \; Y \sim \rho] \leq \theta$ replaces \cref{lem: GUE_norm_bound} from the proof of \cref{thm: no_scaling_shattering}. If the support of $\rho$ is bounded, we can take $\theta = 0$. Otherwise, we need $2n\theta = O(n^{-1})$, which can be accomplished by assuming that $\rho$ has an $n$-dependent scaling, in which case $\|\rho\|_{\infty}$ will also depend on $n$.\footnote{If $\|\rho\|_{\infty}^2/n$ does not decay with $n$, as in the setting where $D^{\rho}$ is the diagonal of a GUE matrix,  we can shrink $\delta$ (i.e., the minimum eigenvalue gap in the proof) by a factor of $\|\rho\|_{\infty}$ to compensate. This will imply shattering with high probability, just for a finer shattering grid and smaller $\epsilon$.}  \\
 \indent In both results, we note an important distinction from the general setting. The proof of \cref{thm: no_scaling_shattering} (or \cref{thm: diagonal_shattering}) is simpler than its  counterpart for general pencils because it does not require a complicated argument to bound eigenvector conditioning; instead, we obtain such a bound directly from the Crawford number and \cref{lem: gamma_eigenvector}. This reflects the fact that the eigenvectors of a definite pencil are, by default, well-conditioned enough for divide-and-conquer, at least as long as $\gamma(A,B)$ is not exponentially small in $n$. 

 \subsubsection{Examples}
\begin{figure}[t]
     \centering
     \includegraphics[width=.95\linewidth]{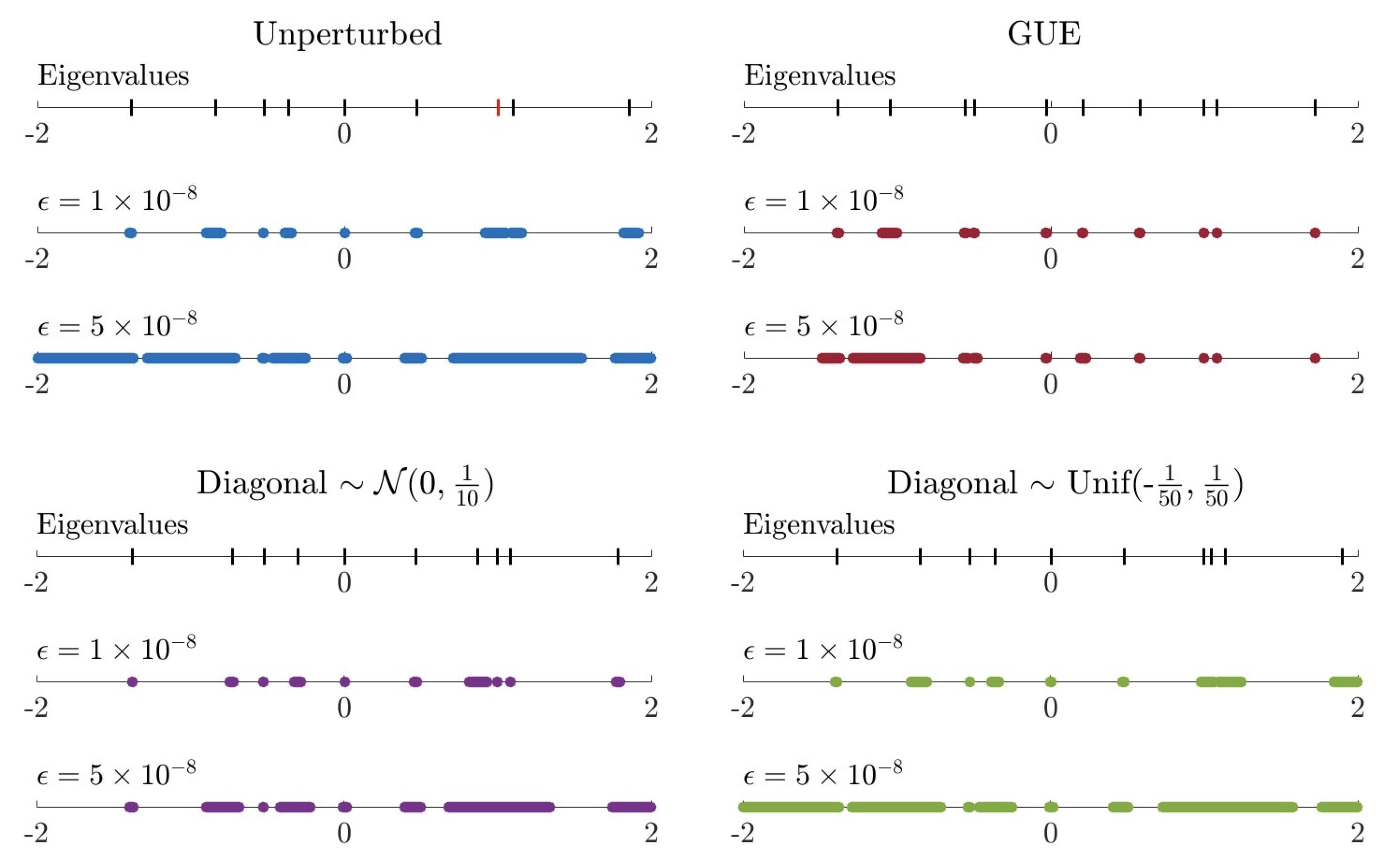}
     \caption{Spectra and symmetric $\epsilon$-pseudospectra of a $10 \times 10$ definite pencil $(A,B)$, constructed as in \eqref{eqn: shattering_examples}, before and after perturbations. Eigenvalues are marked with ticks, while components of $\Lambda_{\epsilon}^{\text{sym}}(A,B)$ are plotted with filled circles. $(A,B)$ initially has a repeated eigenvalue at $+1$, which is highlighted in red, and each perturbation has size $\mu = 10^{-6}$. Initially, $\gamma(A,B) = 5.91 \times 10^{-8}$ while in the remaining subfigures (read clockwise) $\gamma(\widetilde{A},\widetilde{B})$ is $2.44 \times 10^{-7}$, $5.31 \times 10^{-8}$, and $1.04 \times 10^{-7}$, respectively.}
     \label{fig: shattering_comparison}
 \end{figure}

 \indent We provide examples of symmetric pseudospectral shattering under both GUE and diagonal perturbations in \cref{fig: shattering_comparison}. The corresponding definite pencil $(A,B)$ is $10 \times 10$ and constructed from the diagonalization 
 \begin{equation}\label{eqn: shattering_examples}
     (A,B) = (X^H\Lambda X, X^HX),
 \end{equation}
 where $X$ and $\Lambda$ are invertible and diagonal, respectively. To make the pseudospectra of $(A,B)$ more unwieldy, and therefore the shattering effect more dramatic, $X$ and $\Lambda$ are initially constructed randomly and modified as follows:\ a rank-one matrix is subtracted from $X$ to drive up its condition number while two entries of $\Lambda$ are replaced by $+1$ so that $(A,B)$ has a repeated eigenvalue. To plot $\Lambda_{\epsilon}^{\text{sym}}(A,B)$ for this pencil, specifically on the real axis and for values of $\epsilon$ smaller than $\gamma(A,B)$, we leverage \cref{lem: alt_sym_pseudo_equiv}.

 \indent As anticipated, perturbing $(A,B)$ both regularizes its pseudospectra and splits its repeated eigenvalue into two distinct ones, with the size of this effect roughly correlated with the size of the perturbation. For an appropriate set of real grid points, shattering is achieved for $\Lambda_{\epsilon}^{\text{sym}}(\widetilde{A},\widetilde{B})$ with $\epsilon = 10^{-8}$ under both  GUE and diagonal normal perturbations.\footnote{For other combinations of perturbation type and $\epsilon$, shattering may be partially obscured by the relatively large circles used to mark points belonging to the pseudospectra.}    \\
 \indent In all cases, the perturbations applied in \cref{fig: shattering_comparison} are larger than those allowed by \cref{thm: no_scaling_shattering} and \cref{thm: diagonal_shattering}. In particular, we take $\mu = 10^{-6}$ while the algorithm of Kressner, Lu, and Vandereycken \cite{subspace_crawfno} estimates $\gamma(A,B)$ to be $5.82 \times 10^{-8}$. This is problematic in theory only, as the perturbed problems remain definite in all cases. Indeed, two of the perturbations actually increase the Crawford number. This indicates that the restriction on $\mu$ -- which, we recall, is required to guarantee $\gamma(\widetilde{A},\widetilde{B}) \geq \frac{1}{2}\gamma(A,B)$ via standard perturbation theory -- is quite pessimistic. We discuss this further in the next subsection. 
 
\subsubsection{Perturbation Bounds}
Finally, since we will need them in \cref{section: dnc}, we note here a pair of perturbation results for symmetric pseudospectra and symmetric pseudospectral shattering.

\begin{lem}\label{lem: symm_pseudo_under_perturb}
    Let $(A,B)$ be an $n \times n$ definite pencil. If $A',B' \in {\mathbb C}^{n \times n}$ satisfy $\|A - A'\|_2, \|B - B'\|_2 \leq \eta < \frac{\epsilon}{\sqrt{2}}$ for $A - A'$ and $B - B'$ Hermitian, then 
    $$ \Lambda_{\epsilon - \sqrt{2} \eta}^{\text{\normalfont sym}}(A',B') \subseteq \Lambda_{\epsilon}^{\text{\normalfont sym}}(A,B). $$
\end{lem}
\begin{proof}
    Suppose $z \in \Lambda_{\epsilon - \sqrt{2} \eta}^{\text{sym}}(A',B')$. In this case, there exist Hermitian $E,F \in {\mathbb C}^{n \times n}$ with $\|E\|_2^2 + \|F\|_2^2 \leq (\epsilon - \sqrt{2} \eta)^2$ such that $z \in \Lambda(A'+E,B'+F)$. Hence, we have $z \in \Lambda(A + (A'-A+E), B + (B'-B+F))$ with
    \begin{equation}
        \aligned 
        \|A'-A+E\|_2^2 + \|B'-B+F\|_2^2 & \leq (\eta + \|E\|_2)^2 + (\eta + \|F\|_2)^2 \\
        & = 2 \eta^2 + 2 \eta (\|E\|_2 + \|F\|_2) + \|E\|_2^2 + \|F\|_2^2 \\
        & \leq 2 \eta^2 + 2\sqrt{2} \eta (\epsilon - \sqrt{2}\eta) + (\epsilon - \sqrt{2} \eta)^2 \\
        & = \epsilon^2.
        \endaligned 
        \label{eqn: preserve_shattering_one}
    \end{equation}
    Here, the second inequality follows from the fact that $\|E\|_2 + \|F\|_2$ takes maximum value $\sqrt{2}(\epsilon - \sqrt{2}\eta)$ subject to the constraint $\|E\|_2^2 + \|F\|_2^2 \leq (\epsilon - \sqrt{2}\eta)^2$. Since $A'-A+E$ and $B'-B+F$ are Hermitian, we conclude $z \in \Lambda_{\epsilon}^{\text{sym}}(A,B)$. 
\end{proof}

\begin{cor}\label{cor: preserve_definite_shattering}
    Let $(A,B)$ be an $n \times n$ definite pencil and suppose $\Lambda_{\epsilon}^{\text{\normalfont sym}}(A,B)$ is shattered with respect to a set of points $\left\{ g_i \right\} \subset (-r,r)$ for some $0 < \epsilon < \gamma(A,B)$ and $r \in {\mathbb R}_{>0}$. If $A',B' \in {\mathbb C}^{n \times n}$ satisfy 
    $$\|A - A'\|_2 \leq \eta \; \; \; \text{and} \; \; \;  \|B - B'\|_2 \leq \eta$$
    for $\eta < \frac{\epsilon}{\sqrt{2}}$ and $A - A'$ and $B - B'$ Hermitian, then the following hold.
    \begin{enumerate}
        \item Each eigenvalue of $(A',B')$ shares a grid interval $(g_i,g_{i+1})$ with exactly one eigenvalue of $(A,B)$ and $\Lambda_{\epsilon - \sqrt{2}\eta}^{\text{\normalfont sym}}(A',B')$ is also shattered with respect to $g$.
        \item If $(\lambda,v)$ is an eigenpair of $(A,B)$ and 
        $$ \eta  < \frac{\epsilon \gamma(A',B')}{\sqrt{2}\|B\|_2} \cdot \frac{\|B'\|_2 + \|B\|_2}{\|B'\|_2(1 + r^2) + \gamma(A',B')}$$
        then there exists an eigenpair $(\lambda',v')$ of $(A',B')$ such that $\lambda,\lambda' \in (g_i,g_{i+1})$ for some $i$ and 
        $$\frac{\|v - v'\|_2}{\|v\|_2} \leq \frac{\sqrt{2} \eta}{\gamma(A',B')} \cdot \frac{\|B\|_2 \|B'\|_2 (1 + r^2)}{\epsilon \|B'\|_2 + (\epsilon - \sqrt{2}\eta)\|B\|_2} <  1.$$
    \end{enumerate}
\end{cor}
\begin{proof}
    Item (1) follows from a straightforward recreation of the proof of \cite[Lemma 3.14]{arXiv}, noting here that $\gamma(A',B') \geq 1 - \frac{\sqrt{2}\eta}{\gamma(A,B)}$, so $(A',B')$ is definite, and moreover $\Lambda_{\epsilon -\sqrt{2}\eta}^{\text{sym}}(A',B') \subseteq \Lambda_{\epsilon}^{\text{sym}}(A,B)$ by \cref{lem: symm_pseudo_under_perturb}. To prove (2) let $\lambda'$ be the eigenvalue of $(A',B')$ sharing a grid interval with $\lambda$. By construction, any other eigenvalue $\mu \in \Lambda(A',B')$ belongs to a different interval; since $\Lambda_{\epsilon}^{\text{sym}}(A,B)$ and $\Lambda_{\epsilon - \sqrt{2} \eta}^{\text{sym}}(A',B')$ are shattered with respect to $g$ and contain intervals of radius $\frac{\epsilon}{\|B\|_2}$ and $\frac{\epsilon - \sqrt{2}\eta}{\|B'\|_2}$ around the eigenvalues of $(A,B)$ and $(A',B')$, respectively, this guarantees
    \begin{equation}
        |\lambda - \mu| \geq \frac{\epsilon}{\|B\|_2} + \frac{\epsilon - \sqrt{2}\eta}{\|B'\|_2} = \frac{\epsilon \|B'\|_2 + (\epsilon - \sqrt{2}\eta)\|B\|_2}{\|B\|_2\|B'\|_2}.
        \label{eqn: eig_distance}
    \end{equation}
    Consequently, we have
    \begin{equation}
        \chi(\lambda, \mu) = \frac{|\lambda - \mu|}{\sqrt{1 + |\lambda|^2} \sqrt{1 + |\mu|^2}} \geq \frac{\epsilon \|B'\|_2 + (\epsilon - \sqrt{2}\eta) \|B\|_2}{\|B\|_2\|B'\|_2(1+r^2)},
        \label{eqn: bound_chi}
    \end{equation}
    where we note that $|\lambda|, |\mu| \leq r$. Since the criteria on $\eta$ implies 
    \begin{equation}
        \frac{\sqrt{\|A - A'\|_2^2 + \|B - B'\|_2^2}}{\min_{\substack{\mu \in \Lambda(A',B')\\ \mu \neq \lambda'}} \chi(\lambda, \mu)}  \leq \frac{\sqrt{2}\eta \|B\|_2 \|B'\|_2 (1+r^2)}{\epsilon \|B'\|_2 + (\epsilon - \sqrt{2}\eta) \|B\|_2} < \gamma(A',B'),
        \label{eqn: stewart_requirement}
    \end{equation} 
    the result follows from \cref{thm: definite_vector_bound}.
\end{proof}

\subsection{Additional Results for $\gamma(\widetilde{A},\widetilde{B})$}

Before moving on, we pause to consider additional probabilistic bounds for $\gamma(\widetilde{A},\widetilde{B})$. The backdrop here is the complexity discussion from the introduction; we show in the next section that the complexity of our divide-and-conquer routine will depend on a lower bound for the Crawford number of the input pencil, which in our full diagonalization algorithm is $(\widetilde{A},\widetilde{B})$. Hence, we would hope that $\gamma(\widetilde{A},\widetilde{B})$ is, at least, not much smaller than $\gamma(A,B)$, and ideally much larger than the coarse lower bound $\gamma(\widetilde{A},\widetilde{B}) \geq \frac{1}{2} \gamma(A,B)$ used in the proof of shattering. We focus here on the GUE setting (for comments on the diagonal case, see \cref{rem: diag}). \\
\indent We can start by proving an upper bound, which indicates that -- with high probability -- we cannot hope for the perturbation to significantly increase $\gamma(A,B)$ when $n$ is large. Its proof relies on the simple observation that $\frac{1}{\sqrt{2}}(Z_1+iZ_2)$ is a Ginibre random matrix for $Z_1,Z_2 \in \text{GUE}(n)$.

\begin{prop}\label{prop: probabilistic_upper_bound}
    Let $(\widetilde{A},\widetilde{B}) = (A+\mu Z_1, B+\mu Z_2)$ for $(A,B)$ an $n \times n$ Hermitian pencil, $\mu \in {\mathbb R}_{>0}$ and $Z_1,Z_2 \in \text{\normalfont GUE}(n)$. For any $t \geq 0$
    $$ {\mathbb P}\left[\gamma(\widetilde{A},\widetilde{B}) \geq \gamma(A,B)+t \right] \leq e^{-\frac{nt^2}{4\mu^2}}.$$
\end{prop}
\begin{proof}
    Let $y \in {\mathbb C}^n$ be a unit vector satisfying $y^H(A+iB)y = \gamma(A,B)$. We have
    \begin{equation}\label{eqn: prob_lower1}
        \aligned 
        \gamma(\widetilde{A},\widetilde{B}) \leq |y^H(A+iB)y+\mu y^H(Z_1+iZ_2)y| \leq \gamma(A,B) + \mu |y^H(Z_1+iZ_2)y|.
        \endaligned 
    \end{equation}
    Treating $y$ as an arbitrary (fixed) vector, $y^H(Z_1+iZ_2)y$ is normally distributed.\footnote{Because $Z_1$ and $Z_2$ are rotationally invariant, we can assume without loss of generality that $y$ is a standard basis vector, in which case $y^H(Z_1+iZ_2)y$ extracts a diagonal entry of the (scaled) Ginibre matrix $Z_1+iZ_2$.} Hence, we can rewrite \eqref{eqn: prob_lower1} as $\gamma(\widetilde{A},\widetilde{B}) \leq \gamma(A,B) + \mu|Y|$ for a random variable $Y \sim {\mathcal N}_{\mathbb C}(0,\frac{2}{n})$, which implies
    \begin{equation}
        {\mathbb P}[\gamma(\widetilde{A},\widetilde{B}) \geq \gamma(A,B)+t] \leq {\mathbb P}\left[|Y| \geq \frac{t}{\mu}\right] \leq e^{-\frac{nt^2}{4\mu^2}}.
    \end{equation}
    To obtain the final inequality, we note that $|Y|$ follows a chi distribution.
\end{proof}

This upper bound implies the following corollary, which we obtain by setting $A = B = 0$ and $\mu = 1$ in \cref{prop: probabilistic_upper_bound} and integrating. Though not particularly tight, \cref{cor: expected_GUE_craw} captures the fact that the GUE pencil $(Z_1,Z_2)$ is likely to be indefinite for large $n$, which also follows from the circular law \cite{circular}.

\begin{cor}\label{cor: expected_GUE_craw}
    Let $Z_1,Z_2 \in \text{\normalfont GUE}(n)$. ${\mathbb E}[\gamma(Z_1,Z_2)] \leq \sqrt{\pi/n}$. 
\end{cor}

\begin{figure}[t]
    \centering
    \includegraphics[width=\linewidth]{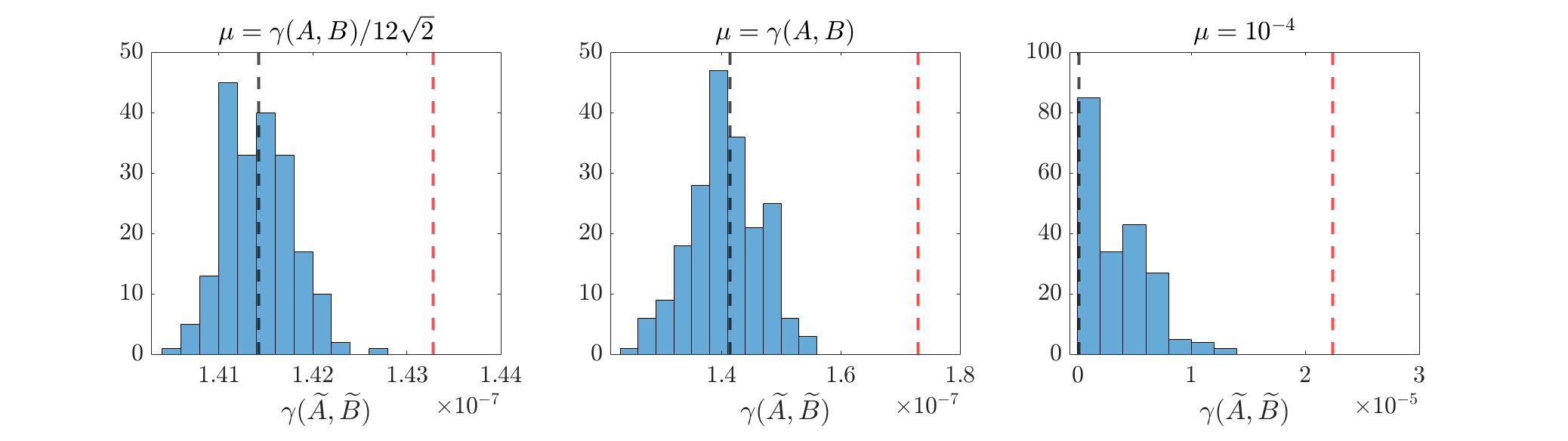}
    \caption{Impact of perturbation size on $\gamma(\widetilde{A},\widetilde{B})$. Each plot presents a histogram of $\gamma(\widetilde{A},\widetilde{B})$ for $200$ GUE perturbations of a $500 \times 500$ definite pencil $(A,B)$, which is constructed so that $\gamma(A,B) = \sqrt{2} \times 10^{-7}$. $\mu$ increases from left to right; at its smallest, it is equal to the bound from \cref{thm: no_scaling_shattering}. We mark $\gamma(A,B)$ and the upper bound $\gamma(A,B)+2\mu\sqrt{\frac{\log(n)}{n}}$ with dashed black and red lines, respectively. The latter holds with probability at least $1-n^{-1}$ by \cref{prop: probabilistic_upper_bound}. In each case, Crawford numbers are computed via the algorithm of Kressner, Lu, and Vandereycken \cite{subspace_crawfno}. }
    \label{fig:crawford_hists}
\end{figure}

\indent Attempting to adapt the proof of \cref{prop: probabilistic_upper_bound} to obtain a lower bound for $\gamma(\widetilde{A},\widetilde{B})$ yields only
\begin{equation}\label{eqn: prob_lower_bound}
    \gamma(\widetilde{A},\widetilde{B}) \geq \gamma(A,B) - \mu \max_{\|x\|_2 = 1} |x^H(Z_1+iZ_2)x|.
\end{equation}
While Collins et al.\ \cite{COLLINS2014516} established that the numerical range of $Z_1+iZ_2$ converges almost surely to the disk of radius two as $n \rightarrow \infty$, implying that the second term of \eqref{eqn: prob_lower_bound} is likely small and that in the limit we have $\gamma(\widetilde{A},\widetilde{B}) \geq \gamma(A,B) - 2\mu$, \eqref{eqn: prob_lower_bound} still yields a lower bound for $\gamma(\widetilde{A},\widetilde{B})$ that is a constant multiple of $\gamma(A,B)$ for the values of $\mu$ allowed by \cref{thm: no_scaling_shattering}. It is for this reason that we stuck with $\gamma(\widetilde{A},\widetilde{B}) \geq \frac{1}{2} \gamma(A,B)$ in our proof of shattering, and for simplicity we will continue to use this bound in the remainder of the paper. \\
\indent Nevertheless, we expect a much tighter result in practice. See for example \cref{fig:crawford_hists}, which records empirical distributions for $\gamma(\widetilde{A},\widetilde{B})$ with varying $\mu$, where $A$ and $B$ are chosen to be simultaneously diagonalizable with a shared eigenvalue at $10^{-7}$ so that $\gamma(A,B) = \sqrt{2} \times 10^{-7}$. We see that for small values of $\mu$, even ones exceeding the restriction of \cref{thm: no_scaling_shattering}, $\gamma(\widetilde{A},\widetilde{B})$ remains close to $\gamma(A,B)$. Indeed, the histograms suggest a rectified-Gaussian-type distribution for $\gamma(\widetilde{A},\widetilde{B})$, though a proof remains out of reach.   \\
\indent At the same time, alternative perturbation strategies yield better lower bounds for $\gamma(\widetilde{A},\widetilde{B})$, at least in some cases. Consider for example random Wishart matrices \cite{wishart}. It is easy to see that the Crawford number of $(A,B)$ can only increase under Wishart perturbations, which are positive definite, if the angle $\theta$ that maximizes \eqref{eqn: alt_crawford_def} satisfies $0 \leq \theta \leq \pi/2$ (equivalently if the point in the field of values of $A+iB$ closest to zero lies in the upper right quadrant -- see the proof of \cite[Theorem VI.1.18]{stewart1990matrix}). As mentioned above, algorithms for computing $\gamma(A,B)$ often provide this information, implying that different perturbation strategies can be used for different input pencils $(A,B)$. Of course, to be viable as a pre-processing step for divide-and-conquer, this would require a proof of shattering in the Wishart case, which remains open.

\begin{remark}\label{rem: diag}
    Much of the preceding analysis can be repeated for diagonal perturbations $D_1^{\rho}$ and $D_2^{\rho}$, though the resulting upper/lower bounds will of course depend on $\rho$. As an example, if $(\widetilde{A},\widetilde{B}) = (A + \mu D_1^{\rho}, B + \mu D_2^{\rho})$ for some (bounded) density $\rho$ on ${\mathbb R}$, then for any $ t > 0$ we have
    \begin{equation}\label{eqn: diag_prob_upper}
        {\mathbb P}\left[ \gamma(\widetilde{A},\widetilde{B}) \geq \gamma(A,B) + t \right] \leq {\mathbb P} \left[ Y^{(n)} \geq \frac{t}{\mu} \right],
    \end{equation}
    where $Y^{(n)}$ is the maximum of $n$ (independent) samples of $|Y_1+iY_2|^2$ for $Y_1$ and $Y_2$ i.i.d.\ according to $\rho$.
\end{remark}

\section{Divide-and-Conquer for Definite Pencils}\label{section: dnc}

We are now ready to construct our specialized divide-and-conquer routine. To build motivation, we start by summarizing the high-level approach for general pencils, as applied in \cite{arXiv}. There, the pencil $(A,B)$ -- which is assumed to come with a guarantee of pseudospectral shattering for the non-symmetrized $\Lambda_{\epsilon}(A,B)$ and a corresponding two-dimensional shattering grid -- is repeatedly split into two smaller pencils as follows:
\begin{enumerate}
    \item First, a grid line is found that separates a fraction of the eigenvalues of $(A,B)$ on either side, which is guaranteed to exist by pseudospectral shattering (recall \cref{fig: shattering_overview}).
    \item Next, projectors onto the deflating subspaces associated to the resulting disjoint sets of eigenvalues are computed. Focusing on one set for simplicity -- e.g., the $k$ eigenvalues above or to the right of the dividing line -- the corresponding right and left projectors $P_R$ and $P_L$ map onto the right and left deflating subspaces ${\mathcal X},{\mathcal Y} \subset {\mathbb C}^n$, where ${\mathcal X}$ is spanned by the associated right eigenvectors of $(A,B)$ and ${\mathcal Y}$ is defined as 
    \begin{equation}
        {\mathcal Y} = \text{span} \left\{ Ax,Bx : x \in {\mathcal X} \right\}.
    \end{equation}
    Informally, ${\mathcal X}$ and ${\mathcal Y}$ generalize the invariant eigenspaces from the standard eigenvalue problem (see \cite{Projector_Paper} for further background).
    \item Finally, computing rank-revealing factorizations of $P_R$ and $P_L$ yields $U_R,U_L \in {\mathbb C}^{n \times k}$, whose orthonormal columns span $\text{range}(P_R)$ and $\text{range}(P_L)$, respectively. The $k \times k$ pencil $(U_L^HAU_R, U_L^HBU_R)$ is then the next subproblem to be divided (in parallel with its counterpart associated to the other set of eigenvalues). Critically, $(\lambda, w)$ is an eigenpair of this subproblem if and only if $(\lambda, U_Rw$) is an eigenpair of $(A,B)$.
\end{enumerate}

\indent Along the way, the $\epsilon$ in our guarantee of pseudospectral shattering provides a benchmark for how accurately $P_R$ and $P_L$ -- and by extension $U_R$ and $U_L$ -- must be computed for shattering to be preserved as we recur. This relies on the fact that $\Lambda_{\epsilon}(U_L^HAU_R,U_L^HBU_R) \subseteq \Lambda_{\epsilon}(A,B)$, assuming $U_L$ and $U_R$ are exact \cite[Lemma 4.1]{arXiv}.

\renewcommand{\arraystretch}{1}
\subsection{Maintaining Definiteness}
We return now to our specialized setting, where the input pencil is definite, with shattering established for $\Lambda_{\epsilon}^{\text{sym}}(A,B)$. If we hope to preserve definiteness through divide-and-conquer, we cannot apply the same splitting procedure outlined above, as $U_L^HAU_R$ and $U_L^HBU_R$ are nonsymmetric in general.\footnote{While the left and right eigenvectors of a definite pencil are the same, left and right deflating subspaces typically are not.} To compensate, we simply replace $U_L$ with another copy of $U_R$. Letting $U = U_R$ to simplify notation, our new subproblem will be $(U^HAU,U^HBU)$, which is clearly symmetric. Moreover, it is definite, with Crawford number at least that of $(A,B)$:
\begin{equation} \label{eqn: subproblem_is_def}
    \gamma(U^HAU, U^HBU) = \min_{\substack{\|y\|_2 = 1 \\ y \in \text{range}(U)}} |y^H(A+iB)y| \geq \gamma(A,B).
\end{equation}
Importantly, the key heuristic of divide-and-conquer also carries over; it is easy to see that the eigenpairs of $(U^HAU, U^HBU)$ are again in correspondence with those of $(A,B)$. \\
\indent Unfortunately, we pay a penalty for dropping $U_L$ in our recursive pseudospectral bound. In particular, we are not guaranteed that $\Lambda_{\epsilon}^{\text{sym}}(U^HAU,U^HBU)$ is contained in $\Lambda_{\epsilon}^{\text{sym}}(A,B)$. Instead, the best we can show is the following.

\begin{lem}\label{lem: definite_pseudo_guarantee}
    Let $(A,B)$ be an $n \times n$ definite pencil and suppose $U \in {\mathbb C}^{n \times k}$ contains an orthonormal basis for a right deflating subspace of $(A,B)$. For any $\epsilon > 0$ and 
    $$ \epsilon' \leq \epsilon \left( 1 + \frac{\|(A,B)\|_2}{\gamma(A,B)} \right)^{-1} $$
    we have $\Lambda_{\epsilon'}^{\text{\normalfont sym}}(U^HAU, U^HBU) \subseteq \Lambda_{\epsilon}^{\text{\normalfont sym}}(A,B)$.
\end{lem}
\begin{proof}
    Suppose $z \in \Lambda_{\epsilon'}^{\text{sym}}(U^HAU,U^HBU)$. In this case, there exist Hermitian matrices $E,F \in {\mathbb C}^{k \times k}$ with $\sqrt{\|E\|_2^2+\|F\|_2^2} \leq \epsilon'$ such that $z \in \Lambda(U^HAU+E, U^HBU+F)$. If $v \in {\mathbb C}^k$ is a corresponding right eigenvector, then by definition
    \begin{equation}
        (U^HAU + E)v = z(U^HBU + F)v.
        \label{eqn: dpg_one}
    \end{equation}
    \indent Consider now $X$, the right eigenvector matrix of $(A,B)$ satisfying $\kappa_2(X) \leq \frac{\|(A,B)\|_2}{\gamma(A,B)}$. Without loss of generality, we may assume that the columns of $U$ span the right deflating subspace corresponding to the first $k$ columns of $X$. Writing $X = [X_1, \; X_2]$ for $X_1 \in {\mathbb C}^{n \times k}$ and $X_2 \in {\mathbb C}^{n \times (n-k)}$, this implies $X_1 = UR_1$ for invertible $R_1 \in {\mathbb C}^{k \times k}$. Storing an orthonormal basis for the range of $X_2$ in $W \in {\mathbb C}^{n \times (n-k)}$, we obtain a block factorization
    \begin{equation}
        X = \begin{pmatrix} U & W \end{pmatrix} \begin{pmatrix} R_1 & 0\\
        0 & R_2 \end{pmatrix}
        \label{eqn: dpg_two}
    \end{equation}
    for another invertible $R_2 \in {\mathbb C}^{(n-k) \times (n-k)}$. \\
    \indent With this in mind, let $Q_R = [U,\; U^{\perp}]$ and $Q_L = [U, \; W]$. Note that $Q_R$ is unitary while $Q_L$ is not (in general). Since the orthogonal complement of the right deflating subspace corresponding to one of $U$ and $W$ is the left deflating subspace associated with the other, we have
    \begin{equation}
        Q_L^HAQ_R = \begin{pmatrix} U^HAU & U^HAU^{\perp} \\ 0& W^HAU^{\perp} \end{pmatrix} \; \;  \text{and} \; \; \; Q_L^HBQ_R = \begin{pmatrix} U^HBU & U^HBU^{\perp} \\ 0& W^HBU^{\perp} \end{pmatrix} .
        \label{eqn: dpg_three}
    \end{equation}
    Hence, it is easy to see 
    \begin{equation}
        \left(Q_L^HAQ_R + \begin{pmatrix} E & 0\\ 0& 0\end{pmatrix} - z \left[ Q_L^HBQ_R + \begin{pmatrix} F & 0\\ 0&0 \end{pmatrix} \right] \right) \begin{pmatrix}  v \\ 0 \end{pmatrix} = 0
        \label{eqn: dpg_four}
    \end{equation}
    or equivalently, noting that $Q_L$ is invertible but not necessarily unitary,
    \begin{equation}
        Q_L^H \left( A + Q_L^{-H} \begin{pmatrix} E & 0 \\ 0& 0 \end{pmatrix} Q_R^H - z \left[ B + Q_L^{-H} \begin{pmatrix} F & 0 \\ 0 & 0 \end{pmatrix} Q_R^H \right] \right) Q_R \begin{pmatrix} v \\ 0 \end{pmatrix} = 0.
        \label{eqn: dpg_five}
    \end{equation}
    In other words, $(z, Uv)$ is an eigenpair of the pencil $(A + M_E, B + M_F)$ for
    \begin{equation}\label{eqn: nonsymmetric_perturb}
        M_E = Q_L^{-H} \begin{pmatrix} E & 0 \\ 0 & 0 \end{pmatrix} Q_R^H \; \; \; \text{and} \; \; \; M_F = Q_L^{-H} \begin{pmatrix} F & 0 \\ 0 & 0 \end{pmatrix} Q_R^H
    \end{equation}
    and therefore
    \begin{equation}\label{eqn: nonsymmetric_equality}
      (A + M_E)Uv =   z (B + M_F)Uv.
    \end{equation}
    \indent While this implies that $z$ belongs to a certain (standard) pseudospectrum of $(A,B)$, it is not sufficient for our purposes since the matrices $M_E$ and $M_F$ are nonsymmetric in general. Instead, consider the following Hermitian perturbations:
    \begin{equation}\label{eqn: hermitian_perturb}
        \aligned 
       \Delta_A &= Q_R\begin{pmatrix} E & U^HM_E^HU^{\perp} \\
        (U^{\perp})^HM_EU & 0 \end{pmatrix} Q_R^H \\
        \Delta_B &= Q_R \begin{pmatrix} F & U^HM_F^HU^{\perp} \\ (U^{\perp})^HM_FU & 0 \end{pmatrix} Q_R^H.
        \endaligned
    \end{equation}
    By construction, applying both \eqref{eqn: dpg_one} and \eqref{eqn: nonsymmetric_equality}, we have
    \begin{equation}\label{eqn: verify_eigenpair}
        \aligned 
        (A + \Delta_A) Uv &= Q_R \left[ Q_R^HAQ_R + \begin{pmatrix} E  & U^H M_E^H U^{\perp} \\ (U^{\perp})^HM_EU & 0 \end{pmatrix} \right] \begin{pmatrix} v \\ 0 \end{pmatrix}  \\
        & = Q_R \begin{pmatrix} U^HAU + E & U^HAU^{\perp} + U^HM_E^HU^{\perp} \\ (U^{\perp})^HAU + (U^{\perp})^HM_EU & (U^{\perp})^HAU^{\perp} \end{pmatrix} \begin{pmatrix} v \\ 0 \end{pmatrix} \\
        & = Q_R \begin{pmatrix} [U^HAU + E]v \\ (U^{\perp})^H [A + M_E]Uv \end{pmatrix} \\
        & = z Q_R \begin{pmatrix}  [U^HBU + F]v \\
        (U^{\perp})^H[B + M_F]Uv \end{pmatrix} \\
        & = z Q_R \begin{pmatrix} U^HBU + F & U^HBU^{\perp} + U^HM_F^HU^{\perp} \\
        (U^{\perp})^HBU + (U^{\perp})^HM_FU & (U^{\perp})^HBU^{\perp} \end{pmatrix} \begin{pmatrix} v \\ 0 \end{pmatrix} \\
        & = z Q_R \left[ Q_R^HBQ_R + \begin{pmatrix} F & U^HM_F^H U^{\perp} \\ (U^{\perp})^HM_FU & 0 \end{pmatrix} \right] \begin{pmatrix} v \\ 0 \end{pmatrix} \\
        & = z(B + \Delta_B)Uv.
        \endaligned 
    \end{equation}
    Thus, $z \in \Lambda(A + \Delta_A, B+\Delta_B)$. \\
    \indent It remains to bound the norms of $\Delta_A$ and  $\Delta_B$. We start by noting that $Q_L^{-1} = \begin{pmatrix} R_1 & 0 \\ 0 & R_2 \end{pmatrix} X^{-1}$ and therefore 
    \begin{equation}
      \|Q_L^{-1}\|_2 \leq \|X^{-1}\|_2 \max \left\{ \|R_1\|_2, \|R_2\|_2 \right\} \leq \kappa_2(X) .
      \label{eqn: dpg_six}
    \end{equation}
    The latter inequality follows from the observation $\|R_1\|_2 = \|X_1\|_2$ and $\|R_2\|_2 = \|X_2\|_2$, which implies that both $\|R_1\|_2$ and $\|R_2\|_2$ are at most $\|X\|_2$. Thus, we have
    \begin{equation} \label{eqn: dpg_seven}
        \|M_E\|_2 = \left| \left| Q_L^{-H} \begin{pmatrix} E & 0 \\ 0 & 0 \end{pmatrix} Q_R^H \right| \right|_2 \leq \kappa_2(X) \|E\|_2  \\
    \end{equation}
    and similarly $\|M_F\|_2 \leq \kappa_2(X)\|F\|_2$. This implies
    \begin{equation}\label{eqn: bound_DeltaA}
        \aligned 
        \|\Delta_A\|_2 &= \left| \left|\begin{pmatrix} E & 0 \\ 0 & 0 \end{pmatrix} + \begin{pmatrix} 0 & U^HM_E^HU^{\perp} \\ (U^{\perp})^HM_EU & 0 \end{pmatrix} \right| \right|_2 \\
        & \leq \|E\|_2 + \|U^HM_E^HU^{\perp}\|_2 \\
        & \leq (1 +  \kappa_2(X)) \|E\|_2
        \endaligned
    \end{equation}\\
    and $\|\Delta_B\|_2 \leq (1 + \kappa_2(X))\|F\|_2 $. Hence, recalling $\kappa_2(X) \leq \frac{\|(A,B)\|_2}{\gamma(A,B)}$ and \cref{defn: symm_pseudospecturum}, to guarantee $z \in \Lambda_{\epsilon}^{\text{sym}}(A,B)$ it is sufficient to take $\epsilon' \leq \epsilon \left( 1 + \frac{\|(A,B)\|_2}{\gamma(A,B)} \right)^{-1}$.
\end{proof}

\cref{lem: definite_pseudo_guarantee} implies that the $\epsilon$ in our guarantee of shattering for $\Lambda_{\epsilon}^{\text{sym}}(A,B)$ must shrink, at least, by a factor of $1+\frac{\|(A,B)\|_2}{\gamma(A,B)}$ at each recursive step. Since here, as in the general setting, the value of $\epsilon$ will determine the amount of work done in the corresponding step of divide-and-conquer, this implies that a structured algorithm may be expensive when $\|(A,B)\|_2/\gamma(A,B)$ is large. In fact, as we will see, \cref{lem: definite_pseudo_guarantee} is the reason the lower bound $\gamma < \gamma(A,B)$ appears in our final complexity (as stated in \cref{section: intro}).\\
\indent Nevertheless, maintaining definiteness has a number of advantages. Most notably, it allows us to assume that the eigenvalues of our input pencil, as well as the eigenvalues of every pencil generated in the divide-and-conquer procedure, are real. Moreover, shattering will be defined relative to a set of grid points instead of grid lines, meaning we have fewer potential splits to check and can guarantee that an optimal one, which divides the eigenvalues into sets of size $\lceil n/2 \rceil$ and $\lfloor n/2 \rfloor$, exists. \\
\indent While these simplifications imply efficiency gains in practice, they cannot on their own reduce the complexity of divide-and-conquer. For example, in both the general and structured settings, we run a binary search to find suitable eigenvalue splits; hence, having roughly half the number of splits to check translates to only a constant factor reduction in the amount of work done. Similarly, splitting each problem closer to 50/50 will reduce the number of recursive steps overall, but the complexity of divide-and-conquer will be the same -- in fact equal to that of its first division -- as long as splits are (at least) any fixed fraction of the current problem size (e.g., no worse than 20/80 as in \cite{arXiv}). \\
\indent To do better, we further leverage the guarantee of real eigenvalues by using a specialized method for computing spectral projectors, swapping the Implicit Repeated Squaring (\textbf{IRS}) routine of \cite{arXiv} (originally introduced by Malyshev \cite{Malyshev1989}) for a recently derived, inverse-free, dynamically weighted Halley iteration (\textbf{IF-DWH}) \cite{Projector_Paper}. The latter computes spectral projectors of $(A,B)$ by implicitly approximating $\text{sign}(B^{-1}A)$. The details of this algorithm (and the matrix sign function) are discussed further in \cref{section: appendix} and \cite{Projector_Paper}. The bottom line is that \textbf{IF-DWH} is provably faster than a routine like \textbf{IRS} when $(A,B)$ has only real eigenvalues. Accordingly, it assumes that $\Lambda(A,B)$ is contained in $[-1,-l) \cup (l,1]$ for an input parameter $l > 0$ (and iteratively drives this union of intervals to $\pm 1$). This kind of bound on the spectrum of $(A,B)$ can be used in tandem with the following result to control approximation error for $\text{sign}(B^{-1}A)$. Recall from \cref{section: background} that $\kappa_2(X)$ can be replaced by $\|(A,B)\|_2/\gamma(A,B)$. For a proof of this result, see \cite[Lemma 3.17]{Projector_Paper}. 

\begin{lem}\label{lem: shifted_dwh_bound}
    Suppose $(A,B)$ is a definite pencil with eigenvalues in $[-1,-l) \cup (l,1]$. If $X$ is any matrix of corresponding right eigenvectors, then
    $$ \|B^{-1}A - \text{\normalfont sign}(B^{-1}A)\|_2 \leq \kappa_2(X) (1-l).$$
\end{lem}

\subsection{Structured Divide-and-Conquer}
We are now ready to state our divide-and-conquer routine, \textbf{EIG-DWH} (\cref{alg: EIG_DWH}) -- to be compared with its general counterpart \textbf{EIG} \cite[Algorithm 5]{arXiv}. 
\begin{algorithm}
\caption[Definite Divide-and-Conquer Eigensolver (\textbf{EIG-DWH})]{Definite Divide-and-Conquer Eigensolver (\textbf{EIG-DWH}) \\
\textbf{Input:} $n \in {\mathbb N}_{+}$, $A,B \in {\mathbb C}^{m \times m}$, $\epsilon, \gamma, c, r  > 0$, $g = \left\{ g_i \right\} \subset {\mathbb R}$ a grid of points, $\eta > 0$ a desired accuracy, and $\theta \in (0,1)$ a failure probability.\\
\textbf{Requires:} $m \leq n$, $(A,B)$ definite with $\gamma(A,B) \geq \gamma$ and $||(A,B)||_2 \leq c$, $g \subset (-r, r)$, and $\Lambda_{\epsilon}^{\text{sym}}(A,B)$ shattered with respect to $g$. \\
\textbf{Output:} $X$ an invertible matrix and $(\Lambda_A,\Lambda_B)$ a diagonal pencil such that $||X^HAX - \Lambda_A||_2 \leq \eta$ and $||X^HBX - \Lambda_B||_2 \leq \eta$ with probability at least $1 - \theta$ (see \cref{thm: EIG_DWH}).}\label{alg: EIG_DWH}
\begin{algorithmic}[1]
\If{$m=1$}
    \State $X = 1$; $\Lambda_A = A$; $\Lambda_B= B$
\Else
    \State $\zeta =  \lfloor \log_2|g|+ 1 \rfloor $
    \State $\delta = \min \left\{ \frac{4\theta}{4 \theta + 3\zeta n^2(n-1)}, \; \sqrt{\frac{\theta}{3(n-1)}} \cdot \frac{\epsilon^2 \gamma^2}{800 n c^2(\gamma+c)^2}, \;\sqrt{\frac{\theta}{3(n-1)}} \cdot \frac{\eta^2}{36nm^2c^2}\right\}$
    \State Choose a grid point $g_i \in g$ 
    \State $({\mathcal A}, {\mathcal B}) = (A - g_iB, 2rB)$
    \State $l = \frac{\epsilon}{2rc}$
    \While{$l < 1 - \frac{2\delta \gamma}{c}$}
        \State $[{\mathcal A},{\mathcal B},l] = \textbf{IF-DWH}({\mathcal A}, {\mathcal B},1,l)$
    \EndWhile
    \State $[U,R_1,R_2, V] = \textbf{GRURV}(2, 2{\mathcal B}, {\mathcal A}+{\mathcal B}, -1, 1)$
    \State $k = \# \left\{ i : \left| \frac{R_2(i,i)}{R_1(i,i)} \right| \geq 2 \sqrt{\frac{\theta}{3 \zeta (n-1)}} \cdot \frac{1-\delta}{n}  \right\}  $ 
    \If{$k < \lfloor \frac{m}{2} \rfloor$ or $k > \lceil \frac{m}{2} \rceil$} 
        \State Return to line 6, executing a binary search over the grid points if necessary.
    \Else 
        \State $U = \textbf{GRURV}(2, 2{\mathcal B}, {\mathcal A}+ {\mathcal B}, -1, 1)$ 
        \State $U_k = U(\; :\; ,1:k)$
        \State $U = \textbf{GRURV}(2, 2{\mathcal B}, {\mathcal A} - {\mathcal B}, -1, 1)$
        \State $U_{m-k} = U(\; : \; , 1:m-k)$
            \State $$ \aligned 
            (A_{11}, B_{11}) &= (U_k^H A U_k, \; U_k^H B U_k ) \\(A_{22}, B_{22}) &= (U_{m-k}^HAU_{m-k}, \; U_{m-k}^HBU_{m-k}) \\
        \endaligned $$
        \State $g_R = \left\{ z \in g : z > g_i \right\}$
        \State $[\widehat{X}, \widehat{\Lambda}_A, \widehat{\Lambda}_B] = \textbf{EIG-DWH}(n, A_{11}, B_{11}, \frac{4\epsilon \gamma}{5(\gamma+c)}, \gamma, c, r,  g_R, \frac{1}{2} \eta, \theta)$
        \State $g_L = \left\{ z \in g : z < g_i \right\}$
        \State $[\widetilde{X}, \widetilde{\Lambda}_A, \widetilde{\Lambda}_B] = \textbf{EIG-DWH}(n, A_{22}, B_{22}, \frac{4\epsilon \gamma}{5(\gamma+c)}, \gamma, c, r, g_L, \frac{1}{2} \eta, \theta)$
        \State $$ 
            X = \left[U_k\widehat{X}, \; U_{m-k}\widetilde{X} \right], \; \; \; \Lambda_A = \text{diag}(\widehat{\Lambda}_A, \widetilde{\Lambda}_A), \; \; \; \Lambda_B = \text{diag}(\widehat{\Lambda}_B,\widetilde{\Lambda}_B)$$
\EndIf
\EndIf
\State \Return $X, \Lambda_A, \Lambda_B$
\end{algorithmic}
\end{algorithm}

\indent We note here the simplifications discussed above:\ \textbf{EIG-DWH} begins by searching over real grid points until a suitable split is found (lines 6-16), computes orthonormal bases for only the corresponding \textit{right} deflating subspaces ($U_k$ and $U_{m-k}$), and subsequently uses them to split the problem in two (while maintaining symmetry). To address the concerns brought up by \cref{lem: definite_pseudo_guarantee}, the $\epsilon$ input parameter shrinks by a factor dependent on $\gamma$ at each step. \\ 
\indent Throughout, the \textbf{GRURV} algorithm of Ballard et al.\ \cite{grurv} is used to (implicitly) compute rank-revealing factorizations. As with \textbf{IF-DWH}, we defer the details of \textbf{GRURV} to \cref{section: appendix}. We do note, however, that \textbf{GRURV} is a randomized algorithm; hence, \textbf{EIG-DWH} is also randomized and takes as input an acceptable failure probability $\theta$.  
\begin{thm}\label{thm: EIG_DWH}
    Let $(A,B)$ and $g = \left\{ g_i \right\} \subset {\mathbb R}$ be a definite pencil and set of grid points satisfying the requirements of \text{\normalfont \textbf{EIG-DWH}} for corresponding parameters $m \leq n$ and $\epsilon, \gamma, c, r > 0$ (see \cref{alg: EIG_DWH}). For any choice of $\eta > 0$ and $\theta \in (0,1)$ suppose
    $$[X,\Lambda_A,\Lambda_B] = \text{\normalfont \bf EIG-DWH}(n,A,B,\epsilon,\gamma,c,r,g,\eta,\theta)$$ 
    in exact arithmetic. With probability at least $1 - \theta$ this call is successful, meaning the recursive procedure converges, and moreover
    $$ \|X^HAX - \Lambda_A\|_2 \leq \eta \; \; \; \; \text{and} \; \; \; \; \|X^HBX - \Lambda_B\|_2 \leq \eta .$$
\end{thm}
\begin{proof}
    This result can be obtained by modifying the proof of \cite[Theorem 4.8]{arXiv}, which demonstrates success for \textbf{EIG}, the general algorithm \textbf{EIG-DWH} is based on. For the sake of brevity, we do not reproduce this proof in full. Instead, we describe the necessary adjustments below.
    \begin{enumerate}
        \item First we note that, as in the general setting, success for \textbf{EIG-DWH} is predicated on the validity of its recursive calls. In this case, there is one additional requirement to check -- i.e., that $(A_{11},B_{11})$ and $(A_{22}, B_{22})$ are definite with $\gamma(A_{11},B_{11}), \gamma(A_{22},B_{22}) \geq \gamma$. This follows from \eqref{eqn: subproblem_is_def}.
        \item Since $\Lambda_{\epsilon}^{\text{sym}}(A,B)$ contains an interval of length $\frac{2\epsilon}{\|B\|_2} \geq \frac{2\epsilon}{c}$ around each eigenvalue (see \cref{thm: symm_BF}), the shattering guarantee $\Lambda_{\epsilon}^{\text{sym}}(A,B) \cap g = \emptyset$ implies that the eigenvalues of $(A-g_iB, B)$ are contained in $ (-2r, -\frac{\epsilon}{c}) \cup (\frac{\epsilon}{c}, 2r)$ for any grid point $g_i \in g$. As a result, each pencil $({\mathcal A}, {\mathcal B})$ input to \textbf{IF-DWH} satisfies $\Lambda({\mathcal A}, {\mathcal B}) \subseteq (-1, -\frac{\epsilon}{2rc}) \cup (\frac{\epsilon}{2rc}, 1)$. This justifies the initial choice of $l$ in line 8.
        \item In line 10, \textbf{IF-DWH} is employed to (implicitly) compute the projector 
        \begin{equation}\label{eqn: grid_point_projector}
            P_{> g_i} = \frac{1}{2}(\text{sign}(B^{-1}A-g_iI)+I).
        \end{equation}
        To guarantee $\|\frac{1}{2} {\mathcal B}^{-1}({\mathcal A}+{\mathcal B}) - P_{> g_i}\|_2 \leq \delta$ upon exit we need $\|{\mathcal B}^{-1}{\mathcal A} - \text{sign}(B^{-1}A-g_iI)\|_2 \leq 2 \delta$. Since
        \begin{equation}
            \text{sign}\left(\frac{1}{2r} B^{-1}(A-g_iB)\right) = \text{sign}(B^{-1}(A-g_iB)) = \text{sign}(B^{-1}A - g_iI),
        \end{equation}
        and moreover noting that \textbf{IF-DWH} changes the eigenvalues of $({\mathcal A}, {\mathcal B})$ but not its eigenvectors or $\text{sign}({\mathcal B}^{-1}{\mathcal A})$, \cref{lem: shifted_dwh_bound} implies that we should run $\textbf{IF-DWH}$ until $l \geq 1 - \frac{2\delta \gamma}{c}$. 
        \item Because the spectrum of $(A,B)$ is confined to the real axis, an optimal split that divides it into disjoint sets of size $\lfloor \frac{m}{2} \rfloor$ and $\lceil \frac{m}{2} \rceil$ always exists. Ensuring optimality at each split guarantees that \textbf{EIG-DWH} is called exactly $n-1$ times on problems of size $m > 1$.
        \item For any $\nu_1 \in (0,1)$ suppose $\delta < \frac{4 \nu_1^2}{4\nu_1^2 + n^2}$ and 
        \begin{equation}
            r = \# \left\{ i : \left| \frac{R_2(i,i)}{R_1(i,i)} \right| \geq \frac{2 \nu_1 (1-\delta)}{n} \right\}
            \label{eqn: r_criteria}
        \end{equation}
        in line 12. In this case, by requiring $k = \text{rank}(P_{> g_i})$ for every grid point checked, a given step of \textbf{EIG-DWH} finds an optimal split with probability at least $1 - \zeta \nu_1^2$, noting that we check at most $\zeta$ grid points (see step four in the proof of \cite[Theorem 3.8]{arXiv}).
        \item Once a split is found, lines 17-20 compute orthonormal bases for the corresponding right deflating subspaces, replacing the calls to \textbf{DEFLATE} \cite[Algorithm 4]{arXiv} in \textbf{EIG}. Accordingly, a straightforward replication of \cite[Theorem 4.7]{arXiv} implies that for $\nu_2 \in (0,1)$ the matrices $U_k$ and $U_{m-k}$ are computed (independently) to within spectral norm error $2 \sqrt{ \frac{n \delta}{\nu_2}}$ with probability at least $1 - 2\nu_2^2$. 
        \item For an appropriate choice of $\delta$, avoiding failure in items 5 and 6 above will guarantee success for one step of divide-and-conquer. Hence, we set $\zeta\nu_1^2 = \nu_2^2 = \frac{\theta}{3(n-1)}$, thereby guaranteeing that a simple union bound will imply a total failure probability for \textbf{EIG-DWH} of at most $\theta$. 
        \item It remains to find the aforementioned choice of $\delta$. Suppose $U \in  {\mathbb C}^{m \times k}$ is the true matrix approximated by $U_k$ -- i.e., satisfying $\|U_k - U\|_2 \leq 2\sqrt{\frac{n \delta}{\nu_2}}$. By \cref{lem: definite_pseudo_guarantee}, $\Lambda_{\epsilon'}(U^HAU, U^HBU)$ is shattered with respect to $g$ for $\epsilon' = \frac{\epsilon \gamma}{\gamma + c}$. Hence, recalling \cref{cor: preserve_definite_shattering}, we maintain shattering for the subsequent calls to \textbf{EIG-DWH} by requiring
        \begin{equation}
            \|A_{11} - U^HAU\|_2, \|B_{11} - U^HBU\|_2 \leq \frac{\epsilon \gamma}{5 \sqrt{2} (\gamma+c)}.
            \label{eqn: def_subproblem_accuracy}
        \end{equation}
        Noting 
        \begin{equation}
            \|A_{11} - U^HAU\|_2 \leq 2 \|U_k-U\|_2 \|A\|_2 \leq 4c \sqrt{ \frac{n \delta}{\nu_2}},
            \label{eqn: def_subproblem_accuracy_2}
        \end{equation}
        it is sufficient to take $4c \sqrt{\frac{n \delta}{\nu_2}} \leq \frac{\epsilon \gamma}{5 \sqrt{2} (\gamma+c)}$ or equivalently $\delta \leq \sqrt{\frac{\theta}{3(n-1)}} \cdot \frac{\epsilon^2 \gamma^2}{800 n c^2(\gamma+c)^2} $, applying our choice of $\nu_2$. The same bounds follow immediately for $(A_{22},B_{22})$.
    \end{enumerate}
    To complete the proof, we show $\|X^HAX - \Lambda_A\|_2 \leq \eta$ and $\|X^HBX - \Lambda_B\|_2 \leq \eta$ under the same conditions that guarantee success for the recursive procedure. This is done inductively. The base case, when $m = 1$, is trivial. For $m > 1$ we have
    \begin{equation}\label{eqn: inductive_diag_one}
        X^HAX - \Lambda_A = \begin{pmatrix}
        \widehat{X}^HA_{11} \widehat{X} - \widehat{\Lambda}_A & \widehat{X}^HU_k^HAU_{m-k}\widetilde{X} \\ \widetilde{X}^HU_{m-k}^HAU_k\widehat{X} & \widetilde{X}^HA_{22}\widetilde{X}-\widetilde{\Lambda}_A        \end{pmatrix}
    \end{equation}
    and therefore
    \begin{equation}\label{eqn: inductive_diag_two}
        \aligned 
        \|X^HAX - \Lambda_A\|_2 \leq \max & \left\{ \|\widehat{X}^HA_{11}\widehat{X} - \widehat{\Lambda}_A\|_2, \; \|\widetilde{X}^HA_{22}\widetilde{X}-\widetilde{\Lambda}_A\|_2 \right\} \\
        & + \|\widehat{X}^HU_k^HAU_{m-k}\widetilde{X}\|_2.
        \endaligned
    \end{equation}
    By our induction hypothesis the first term in this sum is at most $\frac{\eta}{2}$. Since the columns of $\widehat{X}$ and $\widetilde{X}$ have unit length, the remaining term can be bounded as
    \begin{equation}\label{eqn: inductive_diag_three}
        \|\widehat{X}^H U_k^HAU_{m-k}\widetilde{X}\|_2 \leq \|\widehat{X}\|_2\|U_k^HAU_{m-k}\|_2\|\widetilde{X}\|_2 \leq \sqrt{k(m-k)}\|U_k^HAU_{m-k}\|_2.
    \end{equation}
    Now as long as $U_k$ and $U_{m-k}$ are computed accurately enough, $U_k^HAU_{m-k}$ will be nearly zero. To see this, let $U_k$ again approximate $U \in {\mathbb C}^{m \times r}$ and let $U_{m-k}$ approximate $W \in {\mathbb C}^{m \times (m-r)}$. Here, the columns of $U$ and $W$ span disjoint right deflating subspaces of $(A,B)$ with $U_k = U+E$ and $U_{m-k} = W+F$ for some $\|E\|_2,\|F\|_2 \leq \tau$. In this case,
    \begin{equation}\label{eqn: inductive_diag_four}
        U_k^HAU_{m-k} = (U+E)^HA(W+F) = E^HAW + U^HAF+E^HAF
    \end{equation}
    where $U^HAW = 0$ since $(A,B)$ is definite. Consequently,
    \begin{equation}\label{eqn: inductive_diag_five}
        \|U_k^HAU_{m-k}\|_2 \leq (2\tau+\tau^2)\|A\|_2.
    \end{equation}
    If \textbf{EIG-DWH} succeeds as described above we can assume $\tau \leq 2 \sqrt{\frac{n \delta}{\nu_2}}$. Further applying $\|A\|_2 \leq c$ and $\sqrt{k(m-k)} \leq \frac{m}{2}$, this implies $\|U_k^HAU_{m-k}\|_2 \leq 3cm\sqrt{\frac{n \delta}{\nu_2}}$. Hence, the final requirement on $\delta$ guarantees both terms of \eqref{eqn: inductive_diag_two} are bounded by $\frac{\eta}{2}$ and therefore $\|X^HAX - \Lambda_A\|_2 \leq \eta$. We can repeat this argument to obtain $\|X^HBX - \Lambda_B\|_2 \leq \eta$.
\end{proof}

    

We now obtain our final diagonalization algorithm by combining \cref{alg: EIG_DWH} and \cref{thm: EIG_DWH} with the pseudospectral shattering results from \cref{section: shattering}. For simplicity, we use the GUE version of shattering below; nevertheless, \cref{thm: diagonal_shattering} implies that the GUE perturbations $Z_1$ and $Z_2$ in the proof of \cref{cor: fast_diagonalization} can be replaced by $D^{\rho}_1$ and $D^{\rho}_2$ for an appropriate choice of $\rho$.

\begin{cor}\label{cor: fast_diagonalization}
    Let $(A,B)$ be an $n \times n$ definite pencil with $\|A\|_2,\|B\|_2 \leq 1$ and let $\xi < 1$. Assume that a lower bound $\gamma \leq \gamma(A,B)$ is known and moreover $\xi < \frac{n \gamma}{\sqrt{2}}$. Then there exists an exact-arithmetic, randomized algorithm that takes as inputs $(A,B)$ and $\xi$ and outputs invertible $X \in {\mathbb C}^{n \times n}$ and diagonal $\Lambda_A,\Lambda_B \in {\mathbb R}^{n \times n}$ such that 
    $$ \|X^HAX - \Lambda_A\|_2 \leq \xi \; \; \; \; \text{and} \; \; \; \; \|X^HBX - \Lambda_B\|_2 \leq \xi $$
    with probability at least $1 - O(n^{-1})$. This algorithm is inverse-free and requires at most
    $$ O \left( n^{\omega_0}\log\left(n/\xi\right) \log \left( \log \left( n/\xi \right) + \log(n) \log(\gamma^{-1})\right) \right) $$
    operations, where $O(n^{\omega_0})$ is the complexity of $n \times n$ matrix multiplication.
\end{cor}
\begin{proof}
    Consider the following two-step algorithm. First, draw independent $Z_1$ and $Z_2$ from $\text{GUE}(n)$ and set 
    \begin{equation}\label{eqn: step_one_perturb}
        (\widetilde{A},\widetilde{B}) = (A+\mu Z_1,B+\mu Z_2)
    \end{equation}
    for $\mu < \frac{\xi}{12n}$. Next, let $g \subset {\mathbb R}$ be a grid of points constructed as in \cref{thm: no_scaling_shattering} and call
    \begin{equation}\label{eqn: step_two_run_dnc}
        [X,\Lambda_A,\Lambda_B] =  \textstyle \textbf{EIG-DWH}(n,\widetilde{A},\widetilde{B},\frac{\mu^8}{6n^{11}},\frac{\gamma}{2},3,\frac{3n^7+2\mu^6}{2\mu n^5}, g,\frac{\xi}{2},\frac{1}{n}).
    \end{equation}
    We claim that $X,\Lambda_A$, and $\Lambda_B$ satisfy the desired error bounds with probability at least $1 - O(n^{-1})$. \\
    \indent To show this, we first note that the restriction on $\xi$ ensures $\mu < \frac{\gamma(A,B)}{12\sqrt{2}}$. Hence, \cref{thm: no_scaling_shattering} guarantees that with probability at least $1 - O(n^{-1})$ the call to \textbf{EIG-DWH} is valid, meaning the inputs satisfy the listed requirements and in particular $\Lambda_{\epsilon}^{\text{sym}}(\widetilde{A},\widetilde{B})$ is shattered with respect to $g$ for $\epsilon = \frac{\mu^8}{6n^{11}}$.\footnote{Note that the proof of \cref{thm: no_scaling_shattering} also implies the bounds $\gamma(\widetilde{A},\widetilde{B}) \geq \frac{1}{2}\gamma(A,B) \geq \frac{\gamma}{2}$ and $\|(\widetilde{A},\widetilde{B})\|_2 \leq 3$.} When this occurs, \cref{thm: EIG_DWH} further guarantees that with probability at least $1 - n^{-1}$
    \begin{equation}
            \|X^H\widetilde{A}X-\Lambda_A\|_2 \leq \frac{\xi}{2} \; \; \; \; \text{and} \; \; \; \; \|X^H\widetilde{B}X - \Lambda_B\|_2 \leq \frac{\xi}{2}.
    \end{equation}
    By Bayes' theorem, we can assume access to these error bounds via our two-step algorithm with probability at least $1 - O(n^{-1})$. The corresponding bounds for $A$ and $B$ now follow easily; for $A$ in particular we have
    \begin{equation}\label{eqn: A_error_bound}
        \aligned
        \|X^HAX - \Lambda_A\|_2 &= \|X^HAX-X^H\widetilde{A}X+X^H\widetilde{A}X-\Lambda_A\|_2 \\
        & \leq \|X\|_2^2\|A-\widetilde{A}\|_2 + \|X^H\widetilde{A}X - \Lambda_A\|_2 \\
        & \leq 6n\mu + \frac{\xi}{2} \leq \xi.
        \endaligned
    \end{equation}
    Here, we obtain the final inequality by applying both $\|X\|_2 \leq \sqrt{n}$ and  $\|Z_1\|_2 \leq 6$. The former is a consequence of $X$ having columns of unit length while the latter is included in our guarantee of shattering. \\
    \indent To complete the proof we bound the asymptotic complexity of this two-step algorithm. Clearly, shattering can be achieved in $O(n^2)$ operations. The divide-and-conquer procedure, meanwhile, is much more expensive. In the worst case, a single step of divide-and-conquer checks $O(\log(n/\xi))$ grid points as potential eigenvalue splits. If the while loop executed in lines 9-11 of \cref{alg: EIG_DWH} calls \textbf{IF-DWH} at most $p$ times for any grid point checked, the complexity of a single step of divide-and-conquer is at most $O(m^{\omega_0}\log(n/\xi)p)$, assuming access to an $n \times n$ matrix multiplication routine of complexity $O(n^{\omega_0})$. \\
    \indent Now for each grid point checked, the number of calls to $\textbf{IF-DWH}$ is the number of weighted Halley iterations required to map the initial spectral bound $[-1,-\frac{\epsilon}{2rc}) \cup (\frac{\epsilon}{2rc},1]$ to $[-1,-1+\frac{2\delta \gamma}{c}) \cup (1 - \frac{2\delta \gamma}{c},1]$. Following the analysis of \textbf{IF-DWH} in \cite{Projector_Paper}, and noting that $c = O(1)$ while $r = \text{poly}(n,\xi^{-1})$, this is at most 
    \begin{equation}\label{eqn: bound_number_of_Halley_iters}
        O\left( \log \left( \log \left( \frac{n}{\xi \cdot \epsilon } \right) \right) + \log \left( \log \left( \frac{1}{\delta \cdot \gamma} \right) \right) \right).
    \end{equation}
    To simplify, we note that $\delta = \text{poly}(\xi, \epsilon,n^{-1},\gamma)$ so $\eqref{eqn: bound_number_of_Halley_iters}$ becomes $O(\log(\log(\frac{n}{\xi \cdot \epsilon \cdot \gamma})))$. \\
    \indent How small does $\epsilon$ become through the recursion? Let $\epsilon_0 = \frac{\mu^8}{6n^{11}} = \text{poly}(\xi,n^{-1})$ be its initial value. Each step of divide-and-conquer beyond the first sees $\epsilon$ multiplied by $\frac{4\gamma}{5(\gamma+c)}$. Since there are $O(\log(n))$ recursive steps in total we have $\epsilon = (C\gamma)^{\log(n)} \epsilon_0$ at the smallest subproblems, where $C<1$ is a constant. Consequently, $\log(\epsilon^{-1})$ is at most $O(\log(n/\xi) + \log(n)\log(\gamma^{-1}))$ and therefore \eqref{eqn: bound_number_of_Halley_iters} can be bounded by $O(\log(\log(n/\xi)+\log(n)\log(\gamma^{-1})))$. Plugging this in for $p$ above yields the desired complexity bound for the first step of divide-and-conquer. Because each step cuts the problem size in half, it is easy to see that the complexity of the entire divide-and-conquer procedure is asymptotically the same as its first step. 
\end{proof}

As mentioned in the introduction, this result improves on the complexity of divide-and-conquer for arbitrary pencils, as derived in \cite{arXiv}, by replacing a $\log(n/\xi)$ factor with $\log(\log(n/\xi) + \log(n)\log(\gamma^{-1}))$.

\section{Comparison with Existing Methods} \label{section: Examples}
We now present a handful of numerical experiments to test structured divide-and-conquer against existing algorithms for the definite generalized eigenvalue problem. In contrast to the rest of the paper -- which is primarily theoretical and makes the case for divide-and-conquer on the basis of asymptotic complexity -- the results in this section focus on accuracy, where we demonstrate that divide-and-conquer is similarly competitive. All experiments were conducted in Matlab version R2025a. \\
\indent Since the pseudocode of the preceding section is not particularly practical (note, for example, the potentially prohibitively small value of $\delta$ in line 5 of \textbf{EIG-DWH}, the strict stopping criteria for the dynamically weighted Halley iteration, and the unnecessarily demanding choice to search for an optimal split at each level) we implement a simplified version of the algorithm, which admits the modifications listed below. Experiments involving this implementation should be interpreted primarily as a proof of concept. In particular, we note two important caveats. First, by modifying the algorithm we lose access to the theoretical guarantees of \cref{section: dnc}. Indeed, while we set the desired backward diagonalization accuracy to be $\xi = 10^{-15}$ throughout and the divide-and-conquer procedure never failed, we did not meet the corresponding backward diagonalization guarantee of \cref{cor: fast_diagonalization} in any of our tests, though this does not necessarily mean that the resulting eigenvalue approximations are poor. Second, our code is not yet parallelized and therefore cannot compete with highly-optimized methods like QZ in terms of efficiency on small problems. This is our motivation for focusing on accuracy in this section. \\
\indent We now summarize the main adjustments made in our implementation of structured divide-and-conquer (as sketched in the proof of \cref{cor: fast_diagonalization}):
\begin{enumerate}
    \item First, we relax the parameters in the initial call to \textbf{EIG-DWH}, taking $\omega = \epsilon = \mu/n$ for $\mu = \xi/2$, where $\mu$ controls the size of the perturbation as in \eqref{eqn: step_one_perturb}. For simplicity, we consider only GUE perturbations. 
    \item Next, to model the way the algorithm is likely to be used in practice, we only run divide-and-conquer until $m \leq 50$ (as opposed to $m = 1$), at which point we default to QZ.
    \item Since only a few splits will be made, our shattering grid does not need to cover the entire spectrum of $(A,B)$.\footnote{See the discussion in \cite[Section 6.2]{arXiv}.} Accordingly, we focus our search on the potentially much smaller interval $(-\frac{3n^2}{2},\frac{3n^2}{2})$, and we drop the grid construction altogether. Instead, we simply keep track of left/right search endpoints, which are averaged to obtain a candidate split (i.e., $g_i$ in line 6 of \textbf{EIG-DWH}) and updated based on the subsequent rank-revealing factorization. This emulates a binary search over a fixed grid without the memory overhead.
    \item To further simplify the search for a split, we relax the condition on $k$ in line 14 of \textbf{EIG-DWH} to reject one if $k/m < 0.3$ or $k/m > 0.7$.
    \item Because we plan to test the algorithm on pencils with eigenvalues near infinity, we use Implicit Repeated Squaring (\textbf{IRS}) \cite{Projector_Paper} -- run for $\lceil \log_2(m/\epsilon) \rceil$ iterations at each recursive level -- to compute spectral projectors instead of \textbf{IF-DWH}. This is motivated by instability observed in the Halley iteration in \cite{Projector_Paper}. Exploring and/or resolving this issue is left to future work.\footnote{This issue may be fixed by simply using the ``rescaled" version of \textbf{IF-DWH} presented in \cref{section: appendix} (compare with \cite[Algorithm 4]{Projector_Paper}), though a proof remains open.} 
    \item As in \cref{section: shattering}, software of Kressner, Lu, and Vandereycken \cite{subspace_crawfno} is used to estimate $\gamma(A,B)$, which is supplied as input to divide-and-conquer in place of the lower bound $\gamma$.
    \item Finally, since it is unlikely to be necessary in practice, we do not shrink the value of $\epsilon$ recursively. 
\end{enumerate}

\indent In each of the following experiments, we compare structured divide-and-conquer with the following alternatives:
\begin{itemize}
    \item QZ -- the default algorithm for the generalized eigenvalue problem \cite{QZ}, which does not take advantage of definite structure. 
    \item Cholesky + QR -- a two step algorithm that first computes a Cholesky factorization $B = RR^H$, assuming it exists, and subsequently diagonalizes the Hermitian matrix $R^{-1}AR^{-H}$ via reduction to tridiagonal form. 
    \item Cholesky + Jacobi -- a variant of the previous algorithm that handles the second step with Jacobi's method. Our motivation for including this option is work of Davies, Higham, and Tisseur \cite{Cholesky_symm_def}, which demonstrates that it can be much more stable than Cholesky + QR. For background on Jacobi's method, see \cite{CA_Jacobi}.
\end{itemize}
While other algorithms for the definite generalized eigenvalue problem exist in the literature, these are the most accessible. The first two are available in Matlab as \texttt{eig(A,B,`qz')} and \texttt{eig(A,B,`chol')}, respectively. For the third, we use Matlab's built-in Cholesky routine and form $R^{-1}AR^{-H}$ explicitly before calling the implementation of Jacobi's method that accompanies \cite{CA_Jacobi}. \\
\indent In our main experiment, we test these algorithms on a  $1000 \times 1000$ definite pencil constructed as 
\begin{equation}\label{eqn: test_pencil}
    (A,B) = (Q\Lambda_AQ^T,Q\Lambda_BQ^T)
\end{equation}
for $Q$ a random orthogonal matrix and $\Lambda_A$ and $\Lambda_B$ random diagonal. The diagonal entries of $\Lambda_A$ are sampled from ${\mathcal N}_{\mathbb R}(1,1)$ while all but one diagonal entry of $\Lambda_B$ are sampled from $1+|{\mathcal N}_{\mathbb R}(0,1)|$, with the remaining entry set deterministically. We consider 16 trials, in which this deterministic entry is set to $10^{-i}$ for $0 \leq i \leq 15$ (with all other components of \eqref{eqn: test_pencil} fixed). By construction, this varies $\sigma_{\min}(B)$ without altering its other singular values, which allows us to test the sensitivity of each algorithm to the conditioning of $B$. In each trial, $A$ and $B$ are real symmetric\footnote{We restrict to a real problem here since our code for Jacobi's method does not handle the complex Hermitian case.} and simultaneously diagonalizable, with $B$ positive definite. The latter ensures that the aforementioned Cholesky-based algorithms are viable, a requirement that does not apply to divide-and-conquer. Note that the eigenvalues of $(A,B)$ are simply the ratios of the diagonal entries of $\Lambda_A$ and $\Lambda_B$, all but one of which are fixed across the 16 trials. \\
\indent To limit the error introduced in constructing the input matrices $A$ and $B$, we use variable precision arithmetic in Matlab to compute $Q\Lambda_AQ^T$ and $Q\Lambda_BQ^T$ at the start of each trial, which by default keeps 32 significant digits and therefore emulates quad precision. Once computed, $A$ and $B$ are rounded to double precision and passed to the eigensolvers; hence, we can assume that the only error in the input matrices is double-precision roundoff, and we can bound the uncertainty -- in the chordal metric \eqref{eqn: chordal_metrix} -- of the ``exact" eigenvalues obtained by computing the ratios of $\text{diag}(\Lambda_A)$ and $\text{diag}(\Lambda_B)$ symbolically by
\begin{equation}\label{eqn: perturb_bound}
    \texttt{eps} \cdot \frac{\sqrt{\|\texttt{abs}(A)\|_2^2 + \|\texttt{abs}(B)\|_2^2}}{\gamma(A,B)},
\end{equation}
where $\texttt{eps}$ is double-precision machine epsilon and $\texttt{abs}(A)$  is the matrix obtained by taking the absolute value of each entry in $A$. This follows from \cref{thm: definite_value_bound}. 

\begin{figure}[t]
    \centering
    \includegraphics[width=.9\linewidth]{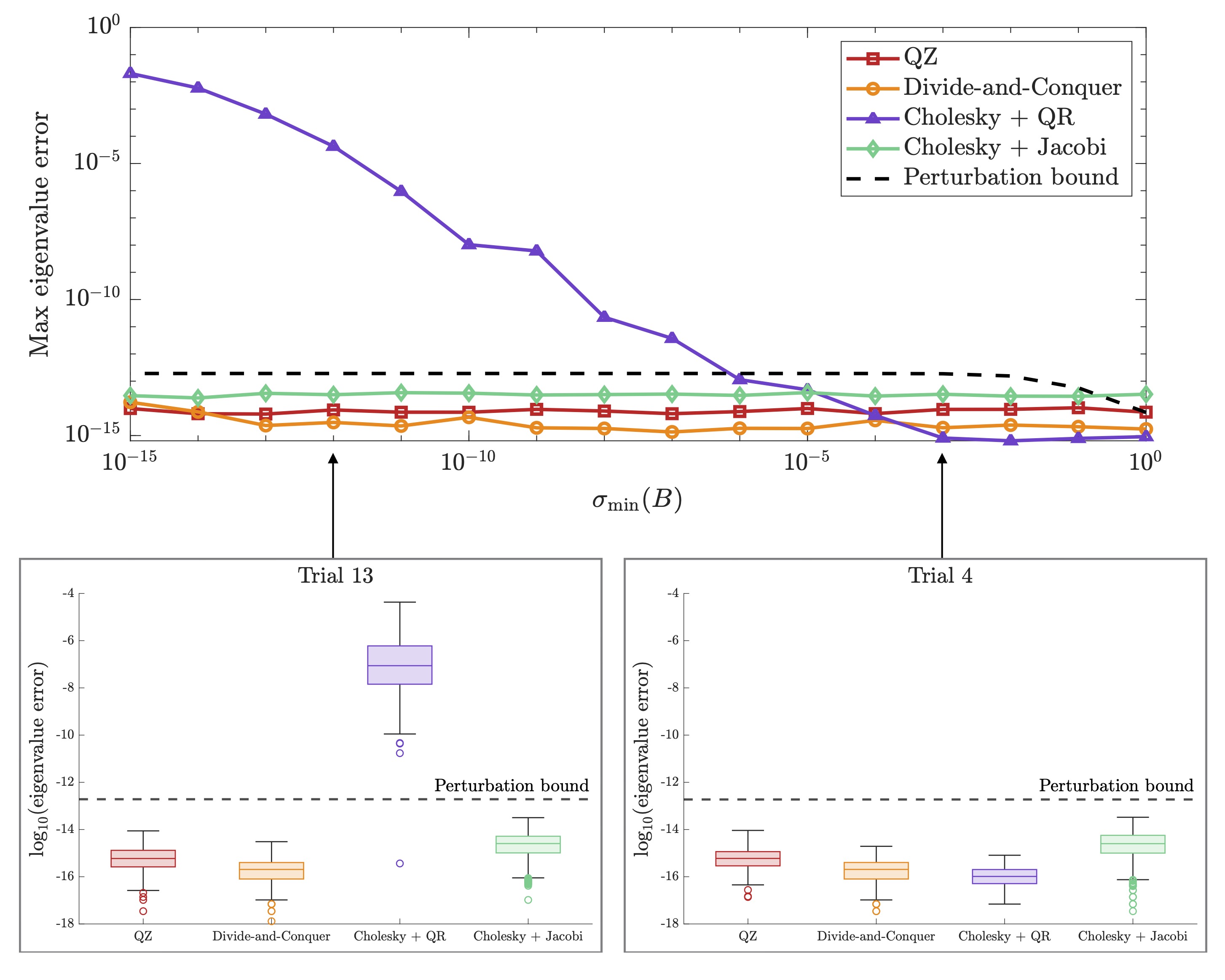}
    \caption{Performance of four algorithms for the definite generalized eigenvalue problem on the $1000 \times 1000$ definite pencil \eqref{eqn: test_pencil}, where we vary the smallest singular value of $B$. Error in each eigenvalue is measured in the chordal metric \eqref{eqn: chordal_metrix}. For context, we also plot the perturbation bound \eqref{eqn: perturb_bound}.}
    \label{fig: main_experiment}
\end{figure}

\indent The results of this experiment are presented in \cref{fig: main_experiment}, which plots the worst-case eigenvalue error, measured in the chordal metric \eqref{eqn: chordal_metrix} relative to the aforementioned ``exact" values, for each algorithm and each trial. We include box plots for two of the trials to confirm that the maximum error accurately captures relative performance across all 1000 eigenvalues. Here, we see that Cholesky + QR is clearly the worst performing algorithm, with accuracy decaying rapidly as $B$ becomes closer to singular. The other three achieve high accuracy regardless of the conditioning of $B$, divide-and-conquer having a slight edge over the other two. \\
\indent In each plot of \cref{fig: main_experiment}, we mark the perturbation bound \eqref{eqn: perturb_bound}. Any method that falls below this dashed curve has computed the eigenvalues as accurately as can be expected in double precision. Note that the shape of the curve indicates that, while initially the Crawford number of $(A,B)$ may shrink with the change in $\sigma_{\min}(B)$, it eventually levels off. In this way, the conditioning of $B$ is decoupled from the definiteness of the pencil in this experiment -- i.e., any instability can be attributed to ill-conditioning in $B$ and not a loss of definiteness. \\
\indent With this in mind, we might ask how the algorithms perform on a pencil made increasingly close to indefinite. We explore this in the next experiment, which runs each method on the $500 \times 500$ pencil
\begin{equation}\label{eqn: test_pencil_two}
    (A,B) = (X\Lambda X^H,X^HX),
\end{equation}
where $X$ is an eigenvector matrix (initially random Gaussian) and $\Lambda$ is random diagonal, with diagonal entries also i.i.d.\ (real) Gaussian. As in the previous experiment, $A$ and $B$ are real symmetric and $B$ is positive definite so that our Cholesky + Jacobi implementation can be used.  

\begin{figure}[t]
    \centering
    \includegraphics[width=.9\linewidth]{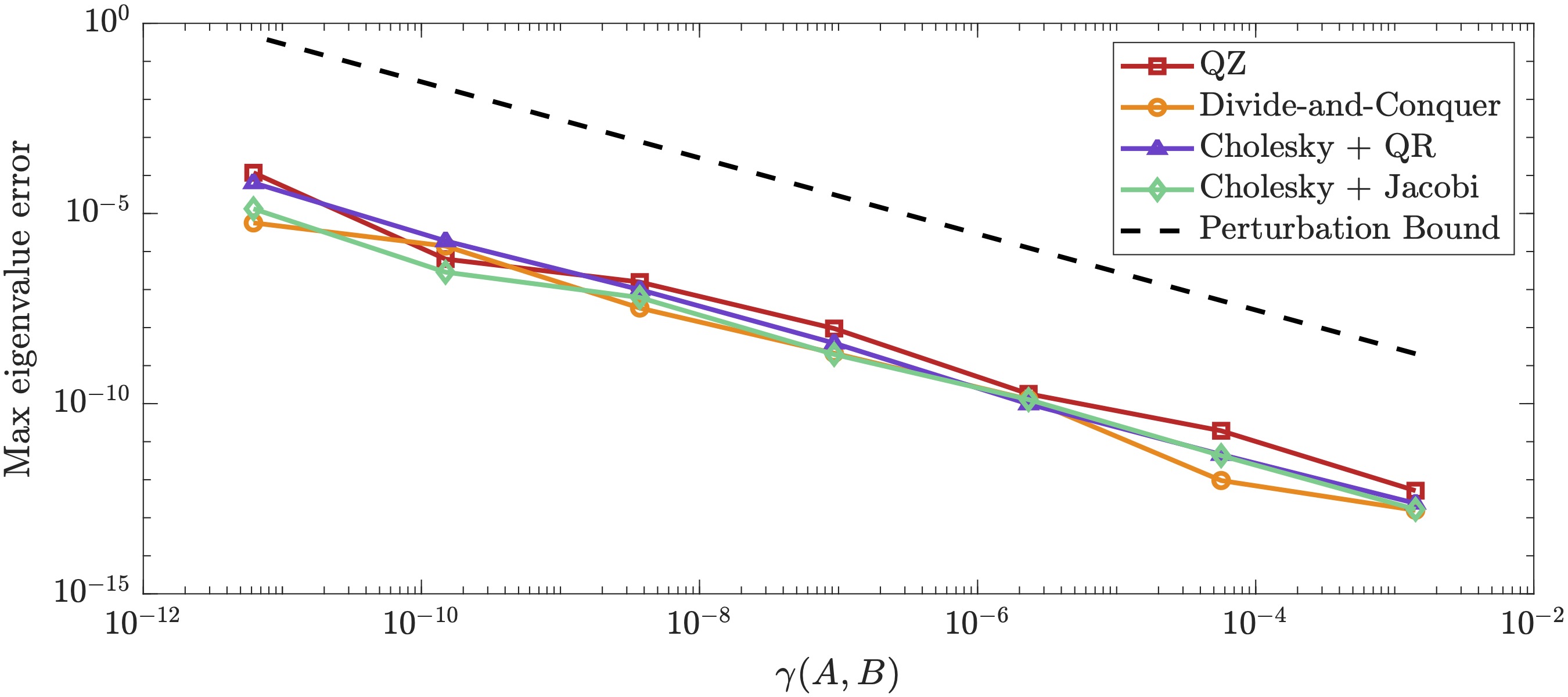}
    \caption{Performance of four algorithms for the definite generalized eigenvalue problem on the $500 \times 500$ pencil \eqref{eqn: test_pencil_two}, which is gradually modified to reduce $\gamma(A,B)$ without changing the corresponding eigenvalues. The perturbation bound again refers to \eqref{eqn: perturb_bound}.}
    \label{fig: experiment_two}
\end{figure}

\indent To reduce the Crawford number of $(A,B)$ over the course of seven trials $1 \leq i \leq 7$, we compute a singular value decomposition $X = U\Sigma V^T$ and update $X$ by subtracting from it the rank one matrix $\frac{5^{i-1}-1}{5^{i-1}} \sigma_{\min}(X)u_nv_n^T$, where $u_n$ and $v_n$ are the last columns of $U$ and $V$, respectively. By construction, this increases the condition number of $X$, which -- for a particular column scaling -- is related to $\gamma(A,B)$ via \cref{lem: gamma_eigenvector}, but does not change the eigenvalues of the pencil. Once again, we use variable precision arithmetic to compute the initial SVD and subsequently to update $X$ and form $A$ and $B$ at the start of each trial. Accordingly, the perturbation bound \eqref{eqn: perturb_bound} again applies to the ``exact" eigenvalues $\text{diag}(\Lambda)$.\\
\indent The results of this experiment are presented in \cref{fig: experiment_two}. As in \cref{fig: main_experiment}, we record the maximum eigenvalue error for each method and each algorithm, as measured in the chordal metric. We see here that all of the methods follow the trend of \eqref{eqn: perturb_bound} -- i.e., they are as sensitive to changes in $\gamma(A,B)$ as we would expect from perturbation theory (and this dominates any loss of accuracy due to ill-conditioning in $B$). \\
\indent Taken together, the results of \cref{fig: main_experiment} and \cref{fig: experiment_two} indicate that structured divide-and-conquer can match (and even beat) the accuracy of existing algorithms for the definite generalized eigenvalue problem. Meanwhile, each of these alternatives comes with one, or more, drawbacks that divide-and-conquer avoids. We note the following:
\begin{enumerate}
    \item Despite being the most widely available algorithm that exploits definiteness, Cholesky + QR is clearly unstable if $B$ is even moderately ill-conditioned.
    \item Cholesky + Jacobi can avoid the instability of Cholesky + QR but is not available if one of $A$ or $B$ is not positive definite. 
    \item While QZ matches the accuracy of divide-and-conquer and can be run on any input pencil, it cannot currently achieve optimal complexity, as established for divide-and-conquer in this paper (arithmetic) and in \cite{Ballard2010MinimizingCF} (communication).
\end{enumerate}

While it may be possible for a better implementation of QZ to narrow the complexity gap with divide-and-conquer, this seems somewhat unlikely. Indeed, there have been several efforts to improve the efficiency of QZ over the years -- e.g., to make better use of level 3 BLAS \cite{QZ2} or in pursuit of parallelism \cite{Parallel_QZ} -- yet none have translated into actual updates to the implementations available in standard libraries like LAPACK,\footnote{Indeed, a LAPACK prospectus still listed this as future work in 2020 \cite{demmel2020prospectus}.} likely due to the difficulty of updating the existing code. While upgrading QZ would no doubt be worthwhile, the results in this paper suggest that implementation efforts would be better spent on divide-and-conquer, which can match QZ in terms of accuracy and is better primed for efficiency on large problems.

\section{Conclusion}\label{section: conclusion}
In this paper, we adapted fast, randomized divide-and-conquer to the definite generalized eigenvalue problem. The result is a highly-parallel, inverse-free algorithm with provably lower complexity than the general version it is based on. \\
\indent For interested readers, we list a handful of open problems and related work:
\begin{enumerate}
    \item As observed in \cref{section: shattering}, different perturbation strategies can be used to prevent the Crawford number from decreasing significantly. Accordingly, it would be useful to derive pseudospectral-shattering-type results for other random (Hermitian) ensembles.
    \item It remains to implement and test \textbf{EIG-DWH} in a truly high-performance fashion. The results of \cref{section: Examples}, together with \cite{banks2020pseudospectral} and \cite{arXiv}, make a compelling case for developing software for randomized divide-and-conquer eigensolvers.
    \item Finally, we note that the definite setting is not the only structured generalized eigenvalue problem worth building a specialized eigensolver for. Most notable here are sparse problems (sparse pseudospectral shattering having been recently established in the single-matrix case by Shah et al.\ \cite{shah2024sparsepseudospectralshattering}). With that in mind, we hope  that this paper provides a roadmap for tailoring divide-and-conquer in other contexts. 
\end{enumerate}

\section{Acknowledgements}\label{section: acknowledgments} 
This work was supported by GFSD and the following NSF grants:\ DMS-2154099, FRG award DMS-1952786, and MSPRF 2402027. Special thanks to Aleksandros Sobczyk, Jorge Garza-Vargas, Frank Uhlig, and Sasha Sodin for helpful correspondence. We also acknowledge feedback from two anonymous referees, which greatly improved the paper.

\bibliographystyle{abbrv}{}
\bibliography{bib}

@article{banks2020pseudospectral,
  title={Pseudospectral {S}hattering, the {S}ign {F}unction, and {D}iagonalization in {N}early {M}atrix {M}ultiplication {T}ime},
  author={Banks, Jess and Garza-Vargas, Jorge and Kulkarni, Archit and Srivastava, Nikhil},
  journal={Foundations of Computational Mathematics},
  pages={1959-2047},
  year={2023},
  volume={23},
  issue={6}, 
  doi={10.1007/s10208-022-09577-5},
}

@article{banks2020gaussian,
author = {Banks, Jess and Kulkarni, Archit and Mukherjee, Satyaki and Srivastava, Nikhil},
year = {2021},
month = {10},
pages = {2114-2131},
title = {Gaussian {R}egularization of the {P}seudospectrum and {D}avies’ {C}onjecture},
volume = {74},
journal = {Communications on Pure and Applied Mathematics},
doi = {10.1002/cpa.22017}
}

@techreport{Ballard2010MinimizingCF,
    Author = {Ballard, Grey and Demmel, James and Dumitriu, Ioana},
    Title = {Minimizing {C}ommunication for {E}igenproblems and the {S}ingular {V}alue {D}ecomposition},
    Institution = {EECS Department, University of California, Berkeley},
    Year = {2011},
    Month = {Feb},
    Number = {UCB/EECS-2011-14},
}

@article{ELSNER1982341,
title = {Perturbation theorems for the generalized eigenvalue problem},
author = {Ludwig Elsner and Sun, Jiguang},
journal = {Linear Algebra and its Applications},
volume = {48},
pages = {341-357},
year = {1982},
issn = {0024-3795},
doi = {10.1016/0024-3795(82)90120-3}
}

@article{2007,
   title={Fast linear algebra is stable},
   volume={108},
   ISSN={0945-3245},
   DOI={10.1007/s00211-007-0114-x},
   journal={Numerische Mathematik},
   publisher={Springer Science and Business Media LLC},
   author={Demmel, James and Dumitriu, Ioana and Holtz, Olga},
   year={2007},
   pages={59–91} 
   }

@article{Bai:CSD-94-793,
  title={An inverse free parallel spectral divide and conquer algorithm for nonsymmetric eigenproblems},
  author={Zhaojun Bai and James Demmel and Ming Gu},
  journal={Numerische Mathematik},
  year={1997},
  volume={76},
  pages={279-308}
}

@article{QZ,
 ISSN = {00361429},
 URL = {http://www.jstor.org/stable/2156353},
 author = {C. B. Moler and G. W. Stewart},
 journal = {SIAM Journal on Numerical Analysis},
 number = {2},
 pages = {241--256},
 publisher = {Society for Industrial and Applied Mathematics},
 title = {An Algorithm for Generalized Matrix Eigenvalue Problems},
 volume = {10},
 year = {1973}
}

@article{grurv,
  doi = {10.48550/ARXIV.1909.06524},
  author = {Ballard, Grey and Demmel, James and Dumitriu, Ioana and Rusciano, Alexander},
  title = {A Generalized Randomized Rank-Revealing Factorization},
  volume = {arXiv:1909.06524},
  year = {2019}
}

@book{TrefethenEmbree+2020,
author = {Lloyd N. Trefethen and Mark Embree},
doi = {10.1515/9780691213101},
title = {Spectra and Pseudospectra: The Behavior of Nonnormal Matrices and Operators},
year = {2020},
publisher = {Princeton University Press},
ISBN = {9780691213101}
}

@article{QZ2,
author = { K\r{a}gstr\"{o}m, Bo and  Kressner, Daniel},
title = {Multishift Variants of the {Q}{Z} Algorithm with Aggressive Early Deflation},
journal = {SIAM Journal on Matrix Analysis and Applications},
volume = {29},
number = {1},
pages = {199-227},
year = {2007},
doi = {10.1137/05064521X}}

@article{Parallel_QZ,
author = {Adlerborn, Bj\"{o}rn and K\r{a}gstr\"{o}m, Bo and Kressner, Daniel},
title = {A Parallel {Q}{Z} Algorithm for Distributed Memory {H}{P}{C} Systems},
journal = {SIAM Journal on Scientific Computing},
volume = {36},
number = {5},
pages = {C480-C503},
year = {2014},
doi = {10.1137/140954817}}

@ARTICLE{SVM,
  author={Mangasarian, O.L. and Wild, E.W.},
  journal={IEEE Transactions on Pattern Analysis and Machine Intelligence}, 
  title={Multisurface proximal support vector machine classification via generalized eigenvalues}, 
  year={2006},
  volume={28},
  number={1},
  pages={69-74},
  doi={10.1109/TPAMI.2006.17}}

@article{MALYSHEV,
title = {Parallel algorithm for solving some spectral problems of linear algebra},
journal = {Linear Algebra and its Applications},
volume = {188-189},
pages = {489-520},
year = {1993},
doi = {10.1016/0024-3795(93)90477-6},
author = {Alexander N. Malyshev}}

@article{BEAVERS1974143,
title = {A new similarity transformation method for eigenvalues and eigenvectors},
author = {A.N. Beavers and E.D. Denman},
journal = {Mathematical Biosciences},
volume = {21},
number = {1},
pages = {143-169},
year = {1974},
issn = {0025-5564},
doi = {10.1016/0025-5564(74)90111-4}}

@article{circular,
author = {Z. D. Bai},
title = {{Circular law}},
volume = {25},
journal = {The Annals of Probability},
number = {1},
publisher = {Institute of Mathematical Statistics},
pages = {494-529},
year = {1997},
doi = {10.1214/aop/1024404298}}

@book{stewart1990matrix,
  title={Matrix Perturbation Theory},
  author={Stewart, G. W. and Sun, J.},
  isbn={9780126702309},
  lccn={lc90033378},
  series={Computer Science and Scientific Computing},
  year={1990},
  publisher={Elsevier Science}
}

@article{Aizenman17,
author = {Aizenman, Michael and Peled, Ron and Schenker, Jeffrey and Shamis, Mira and Sodin, Sasha},
title = {Matrix regularizing effects of {G}aussian perturbations},
journal = {Communications in Contemporary Mathematics},
volume = {19},
number = {03},
pages = {1750028},
year = {2017},
doi = {10.1142/S0219199717500286}}

@book{Higham,
author = {Higham, Nicholas J.},
title = {Accuracy and Stability of Numerical Algorithms},
publisher = {Society for Industrial and Applied Mathematics},
year = {2002},
doi = {10.1137/1.9780898718027},
edition = {{S}econd},
}

@article{arXiv,
author = {Demmel, James and Dumitriu, Ioana and Schneider, Ryan},
year = {2024},
title = {Generalized Pseudospectral Shattering and Inverse-Free Matrix Pencil Diagonalization},
journal = {Foundations of Computational Mathematics},
doi = {10.1007/s10208-024-09682-7},
}

@article{bulgakov,
title={Circular dichotomy of the spectrum of a matrix},
author={A. Ya. Bulgakov and S.K. Godunov},
year={1988},
journal={Siberian Mathematical Journal},
volume={29},
pages={734-744},
doi={10.1007/BF00970267}
}

@article{Malyshev1989,
author = {A. N. Malyshev},
year = {1989},
title = {Computing invariant subspaces of a regular linear pencil of matrices},
journal =  {Siberian Mathematical Journal},
pages = {559--567},
volume = {30},
doi= {10.1007/BF00971756}
}

@article{Projector_Paper,
title = {Fast and inverse-free algorithms for deflating subspaces},
journal = {Linear Algebra and its Applications},
volume = {735},
pages = {222-258},
year = {2026},
issn = {0024-3795},
doi = {10.1016/j.laa.2026.01.014},
author = {James Demmel and Ioana Dumitriu and Ryan Schneider},
}

@article{shah2024sparsepseudospectralshattering,
      title={Sparse Pseudospectral Shattering}, 
      author={Rikhav Shah and Nikhil Srivastava and Edward Zeng},
      year={2024},
      volume={arXiv:2411.19926},
      doi = {10.48550/arXiv.2411.19926}
}

@article{STEWART_DEFINITE,
title = {Perturbation bounds for the definite generalized eigenvalue problem},
journal = {Linear Algebra and its Applications},
volume = {23},
pages = {69-85},
year = {1979},
issn = {0024-3795},
doi = {10.1016/0024-3795(79)90094-6},
author = {Stewart, G. W.},
}

@article{optimizing_halley,
author = {Nakatsukasa, Yuji and Bai, Zhaojun and Gygi, Fran\c{c}ois},
title = {Optimizing {H}alley's {I}teration for {C}omputing the {M}atrix {P}olar {D}ecomposition},
journal = {SIAM Journal on Matrix Analysis and Applications},
volume = {31},
number = {5},
pages = {2700-2720},
year = {2010},
doi = {10.1137/090774999},
}

@article{wishart,
 author = {John Wishart},
 journal = {Biometrika},
 number = {1/2},
 pages = {32--52},
 publisher = {[Oxford University Press, Biometrika Trust]},
 title = {The Generalized Product Moment Distribution in Samples from a Normal Multivariate Population},
 volume = {20A},
 year = {1928},
 doi = {10.2307/2331939}
}

@article{Minami,
author = {Nariyuki Minami},
title = {{Local fluctuation of the spectrum of a multidimensional Anderson tight binding model}},
volume = {177},
journal = {Communications in Mathematical Physics},
number = {3},
publisher = {Springer},
pages = {709-725},
year = {1996},
doi = {10.1007/BF02099544},
}

@article{Crawford_1,
 author = {C. R. Crawford},
 journal = {SIAM Journal on Numerical Analysis},
 number = {6},
 pages = {854--860},
 publisher = {Society for Industrial and Applied Mathematics},
 title = {A {S}table {G}eneralized {E}igenvalue {P}roblem},
 volume = {13},
 year = {1976}
}

@article{Crawford_thesis,
author={Crawford, Charles R.},
year={1970},
title={The {N}umerical {S}olution of the {G}eneralized {E}igenvalue {P}roblem},
volume={PhD Thesis},
}

@article{My_thesis,
author={Schneider, Ryan},
year={2024},
title={Pseudospectral {D}ivide-and-{C}onquer for the {G}eneralized {E}igenvalue {P}roblem},
volume={PhD Thesis},
}

@article{Crawford_Moon,
title = {Finding a positive definite linear combination of two {H}ermitian matrices},
journal = {Linear Algebra and its Applications},
volume = {51},
pages = {37-48},
year = {1983},
doi = {10.1016/0024-3795(83)90148-9},
author = {C. R. Crawford and Yiu Sang Moon},
}

@article{arg_alg,
author = {Guo, Chun-Hua and Higham, Nicholas J. and Tisseur, Fran\c{c}oise},
title = {An {I}mproved {A}rc {A}lgorithm for {D}etecting {D}efinite {H}ermitian {P}airs},
journal = {SIAM Journal on Matrix Analysis and Applications},
volume = {31},
number = {3},
pages = {1131-1151},
year = {2010},
doi = {10.1137/08074218X},
}

@article{Uhlig,
title = {On computing the generalized {C}rawford number of a matrix},
journal = {Linear Algebra and its Applications},
volume = {438},
number = {4},
pages = {1923-1935},
year = {2013},
doi = {10.1016/j.laa.2011.06.024},
author = {Frank Uhlig},
}

@article{FORD1974337,
title = {The generalized eigenvalue problem in quantum chemistry},
journal = {Computer Physics Communications},
volume = {8},
number = {5},
pages = {337-348},
year = {1974},
issn = {0010-4655},
doi = {10.1016/0010-4655(74)90011-3},
author = {Brian Ford and George Hall},
}

@article{Wegner,
author = {Wegner, Franz},
year = {1981},
title = {Bounds on the density of states in disordered systems},
journal = {Zeitschrift für Physik B Condensed Matter},
pages = {9-15},
volume = {44},
issue = {1},
doi = {10.1007/BF01292646}
}

@article{COLLINS2014516,
title = {Numerical range for random matrices},
journal = {Journal of Mathematical Analysis and Applications},
volume = {418},
number = {1},
pages = {516-533},
year = {2014},
doi = {10.1016/j.jmaa.2014.03.072},
author = {Benoît Collins and Piotr Gawron and Alexander E. Litvak and Karol Życzkowski},
}

@article{Veselié1993,
author = {Veselić, K},
year = {1993},
title = {A {J}acobi eigenreduction algorithm for definite matrix pairs},
journal = {Numerische Mathematik},
pages = {241-269},
volume = {64},
issue = {1},
doi = {10.1007/BF01388689}
}

@article{Mehl_jacobi,
author = {Mehl, Christian},
title = {Jacobi-like Algorithms for the Indefinite Generalized {H}ermitian Eigenvalue Problem},
journal = {SIAM Journal on Matrix Analysis and Applications},
volume = {25},
number = {4},
pages = {964-985},
year = {2004},
doi = {10.1137/S089547980240947X}, 
}

@article{Chandrasekaran00,
author = {Chandrasekaran, S.},
title = {An Efficient and Stable Algorithm for the Symmetric-Definite Generalized Eigenvalue Problem},
journal = {SIAM Journal on Matrix Analysis and Applications},
volume = {21},
number = {4},
pages = {1202-1228},
year = {2000},
doi = {10.1137/S0895479897316308},
}

@article{Bunse84,
author = {A. Bunse-Gerstner},
title = {An algorithm for the symmetric generalized eigenvalue problem},
journal = {Linear Algebra and its Applications},
volume = {58},
year = {1984},
page = {43-68},
doi = {10.1016/0024-3795(84)90203-9}
}

@article{Cholesky_symm_def,
author = {Davies, Philip I. and Higham, Nicholas J. and Tisseur, Fran\c{c}oise},
title = {Analysis of the {C}holesky Method with Iterative Refinement for Solving the Symmetric Definite Generalized Eigenproblem},
journal = {SIAM Journal on Matrix Analysis and Applications},
volume = {23},
number = {2},
pages = {472-493},
year = {2001},
doi = {10.1137/S0895479800373498},
}

@article{Scott81,
author = {Scott, David S.},
title = {Solving Sparse Symmetric Generalized Eigenvalue Problems without Factorization},
journal = {SIAM Journal on Numerical Analysis},
volume = {18},
number = {1},
pages = {102-110},
year = {1981},
doi = {10.1137/0718008},
}

@article{subspace_crawfno,
author = {Kressner, Daniel and Lu, Ding and Vandereycken, Bart},
title = {Subspace Acceleration for the {C}rawford Number and Related Eigenvalue Optimization Problems},
journal = {SIAM Journal on Matrix Analysis and Applications},
volume = {39},
number = {2},
pages = {961-982},
year = {2018},
doi = {10.1137/17M1127545},
}

@article{GRAILLAT200668,
title = {A note on structured pseudospectra},
journal = {Journal of Computational and Applied Mathematics},
volume = {191},
number = {1},
pages = {68-76},
year = {2006},
issn = {0377-0427},
doi = {10.1016/j.cam.2005.04.027},
author = {Stef Graillat},
}

@article{CA_Jacobi,
      title={Minimizing the Arithmetic and Communication Complexity of {J}acobi's Method for Eigenvalues and Singular Values: Part One - Serial Algorithms}, 
      author={James Demmel and Hengrui Luo and Ryan Schneider and Yifu Wang},
      year={2025},
      volume={arXiv:2506.03466},
      doi={10.48550/arXiv.2506.03466},
}

@article{structured_dnc_arxiv,
      title={Structured Divide-and-Conquer for the Definite Generalized Eigenvalue Problem}, 
      author={James Demmel and Ioana Dumitriu and Ryan Schneider},
      year={2025},
      volume={arXiv:2505.21917},
      doi = {10.48550/arXiv.2505.21917},
}

@article{demmel2020prospectus,
  title={Prospectus for the Next {LAPACK} and {S}ca{LAPACK} Libraries: Basic ALgebra Libraries for Sustainable Technology with Interdisciplinary Collaboration ({BALLISTIC})},
  author={Demmel, James and Dongarra, Jack and Langou, Julie and Langou, Julien and Luszczek, Piotr and Mahoney, Michael W},
  journal={University of Tennessee, Knoxville, TN},
  year={2020}
}

\appendix

\section{GRURV and IF-DWH}\label{section: appendix}

This appendix discusses the numerical details of two building blocks of our structured divide-and-conquer algorithm:\ \textbf{GRURV}, which computes rank-revealing factorizations of products of matrices and their inverses, and \textbf{IF-DWH}, which iteratively approximates $\text{sign}(B^{-1}A)$. Note that both algorithms work implicitly and are therefore inverse free. \\
\indent We start with \textbf{GRURV} (\cref{alg:GRURV}). Originally introduced by Ballard et al.\ \cite{grurv} and built on top of the earlier \textbf{RURV}/\textbf{RULV} algorithm of Demmel, Dumitriu, and Holtz \cite{2007} (presented here as \cref{alg:RURV}), \textbf{GRURV} computes a rank-revealing factorization of $A_1^{m_1}A_2^{m_2}\cdots A_k^{m_k}$ for arbitrary matrices $A_1, \ldots, A_k$ and exponents $m_i$, which are $\pm 1$. As noted in \cref{section: dnc}, \textbf{GRURV} is randomized, built around the fact that, with high probability, a simple QR factorization is rank-revealing when applied to a matrix hit by a random Haar transformation. \textbf{GRURV} is also backwards stable and capable of producing strongly rank-revealing factorizations (see the discussion in \cite{grurv}).\\
\indent In the context of divide-and-conquer, \textbf{GRURV} is particularly compatible with methods that approximate spectral projectors implicitly. This includes \textbf{IF-DWH} (\cref{alg:IF_DWH}). The heuristic behind the approach employed by \textbf{EIG-DWH} is simple:\ $\frac{1}{2} (\text{sign}(A)+I)$ is a projector onto the eigenvectors of $A$ corresponding to eigenvalues in the right half plane, where if $A = P \begin{pmatrix} J_+ & \\& J_- \end{pmatrix}P^{-1}$ is the Jordan canonical form of $A$ (with $J_+$ and $J_-$ associated to eigenvalues in the right and left half planes, respectively) then
\begin{equation}\label{eqn: sign}
    \text{sign}(A) = P \begin{pmatrix} I & \\& -I \end{pmatrix} P^{-1}
\end{equation}
for the same block structure. Assuming $B$ is invertible, $\frac{1}{2}(\text{sign}(B^{-1}A)+I)$ is therefore a projector onto a right deflating subspace of the pencil $(A,B)$, again corresponding to eigenvalues in the right half plane. 

\begin{algorithm}[t]
\caption{Randomized Rank-Revealing Factorization (\textbf{RURV}/\textbf{RULV})\\
\textbf{Input:} $A \in {\mathbb C}^{n \times n}$ \\
\textbf{Output:} $U$ unitary, $R$ upper triangular or $L$ lower triangular, and $V$ Haar such that $A = URV$ or $A = ULV$ is a rank-revealing factorization of $A$. }\label{alg:RURV}
\begin{algorithmic}[1]
\State Draw a random matrix $B$ with i.i.d. ${\mathcal N}_{\mathbb C}(0,1)$ entries
\State $[V, \hat{R}] = \textbf{QR}(B)$
\State $\hat{A} = A \cdot V^H$
\State $[U,R] = \textbf{QR}(\hat{A}) $ or $[U,L] = \textbf{QL}(\hat{A})$
\State \Return $U,R$ or $L$, and $V$
\end{algorithmic}
\end{algorithm}

\begin{algorithm}[t]
\caption{Generalized Randomized Rank-Revealing Factorization (\textbf{GRURV})\\
\textbf{Input:} $k$ a positive integer, $A_1, A_2, \ldots, A_k \in {\mathbb C}^{n \times n}$, and $m_1, m_2, \ldots, m_k \in \left\{ 1, -1 \right\} $ \\
\textbf{Output:} $U$ unitary, $R_1, R_2, \ldots, R_k$ upper triangular, and $V$ Haar such that  $UR_1^{m_1} R_2^{m_2} \cdots R_k^{m_k}V$ is a rank-revealing factorization of $A_1^{m_1}A_2^{m_2} \cdots A_k^{m_k}$ }\label{alg:GRURV}
\begin{algorithmic}[1]
\If{$m_k = 1$}
    \State $[U,R_k,V] = \textbf{RURV}(A_k)$
\Else 
    \State $[U, L_k, V] = \textbf{RULV}(A_k^H)$
    \State $R_k = L_k^H$
\EndIf
\State $U_{\text{current}} = U$
\For{$i = k-1: 1$} 
    \If{$m_i = 1$} 
        \State $[U, R_i] = \textbf{QR}(A_i \cdot U_{\text{current}}) $
        \State $U_{\text{current}} = U$
    \Else
        \State $[U, R_i] = \textbf{RQ}(U_{\text{current}}^H \cdot A_i)$
        \State $U_{\text{current}} = U^H$
    \EndIf
\EndFor
\State \Return $U_{\text{current}}$, optionally $R_1, R_2, \ldots, R_k$, $V$
\end{algorithmic}
\end{algorithm}

\indent In computing projectors this way, \textbf{IF-DWH} is used to approximate $\text{sign}(B^{-1}A)$. Roughly speaking, \textbf{IF-DWH} evaluates (implicitly) a rational function at $B^{-1}A$, thereby producing a new pencil $(A_p,B_p)$ whose right eigenvectors are the same as $(A,B)$ but whose eigenvalues have been driven to $\pm 1$, depending on their relation to the imaginary axis. Efficiency is derived from the dynamic coefficients computed in lines 5-7, which speed up this process for real eigenvalues and were originally found by Nakatsukasa, Bai, and Gygi \cite{optimizing_halley}. Recalling \cref{lem: shifted_dwh_bound}, this implies accurate approximations of $\text{sign}(B^{-1}A)$ after a modest number of iterations. In fact, the convergence properties of \textbf{IF-DWH} are provably better (when applied to pencils with real eigenvalues) than other inverse-free algorithms for spectral projectors (see \cite{Projector_Paper} for the details).\\
\indent In \textbf{EIG-DWH}, these subroutines are used in tandem as follows:
\begin{enumerate}
    \item After a grid point is selected, the pencil $(A,B)$ is shifted and scaled to obtain $({\mathcal A}, {\mathcal B}) = (A-g_iB,2rB)$, whose eigenvalues live in a union of intervals $(-1,-l) \cup (l,1)$. The eigenvalues of $(A,B)$ to the right of $g_i$ correspond to eigenvalues of $({\mathcal A},{\mathcal B})$ in the right half plane. 
    \item Applying \textbf{IF-DWH} to $({\mathcal A},{\mathcal B})$ produces a pencil $(A_p,B_p)$ that approximates $\text{sign}({\mathcal B}^{-1}{\mathcal A})$, implicitly, as $B_p^{-1}A_p$. Since the (right) eigenvectors of $({\mathcal A},{\mathcal B})$ and $(A,B)$ are the same, this simultaneously approximates the spectral projector of $(A,B)$ corresponding to eigenvalues to the right of $g_i$. 
    \item Finally, \textbf{GRURV} computes a rank-revealing factorization of this approximation, which can be expressed as $\frac{1}{2}B_p^{-1}(A_p+B_p)$. Reading the rank from this factorization (which can be done without explicitly inverting the triangular matrices produced by \textbf{GRURV}) allows us to check the quality of $g_i$ as a splitting point before using it to divide the problem.
    \item If the split is selected, a projector onto the right deflating subspace corresponding to the remaining eigenvalues can be obtained without calling \textbf{IF-DWH} a second time. This is a simple consequence of the fact that $\frac{1}{2} (\text{sign}(A)-I)$ is a projector onto eigenvectors of $A$ associated to the left half plane. 
\end{enumerate}

\indent Finally, note that the call to \textbf{GRURV} in line 17 of  \textbf{EIG-DWH} is technically not necessary. We could instead use the $U$ matrix output in line 12. Nevertheless, we keep it in to simplify the argument used to bound the failure probability of the algorithm -- i.e., it allows us to assume that the accuracy of $U$ and $k$ are independent.

\renewcommand{\arraystretch}{1}
\begin{algorithm}[t]
\caption[Inverse-Free Dynamically Weighted Halley Iteration (\textbf{IF-DWH})]{Inverse-Free Dynamically Weighted Halley Iteration (\textbf{IF-DWH})\\
\textbf{Input:} $A, B \in {\mathbb C}^{n \times n}$, $p$ a number of iterations, $l_0 > 0$. \\
\textbf{Requires:} Eigenvalues $\lambda$ of $(A,B)$ are real with $ l_0 < |\lambda|  \leq 1$.}\label{alg:IF_DWH}
\begin{algorithmic}[1]
\State $A_0 = A$
\State $B_0 = B$
\For{$j = 0:p-1$} 
    \State $\gamma_j = \left(4 (1 - l_j^2)/l_j^4 \right)^{1/3} $
    \State $b_j = \sqrt{1 + \gamma_j} + \frac{1}{2} \sqrt{8 - 4 \gamma_j + 8(2 - l_j^2)/(l_j^2 \sqrt{1 + \gamma_j})}$
    \State $a_j = \frac{1}{4} (b_j - 1)^2$
    \State $c_j = a_j + b_j - 1$
    \vspace{1mm}
    \State $\begin{pmatrix} -B_j \\ \sqrt{c_j}A_j
    \end{pmatrix}  = \begin{pmatrix} Q_{11} & Q_{12} \\
    Q_{21} & Q_{22} \end{pmatrix} \begin{pmatrix} R_j \\ 0  \end{pmatrix}$ \Comment{Apply Halley iteration}
    \vspace{1mm}
    \State $C_j = \frac{a_j}{\sqrt{c_j}} \cdot Q_{12}^HA_j + b_j \cdot Q_{22}^HB_j$
    \State $D_j = \sqrt{c_j} \cdot Q_{12}^HA_j + Q_{22}^HB_j$
    \vspace{1mm}
    \State $\begin{pmatrix} -D_j \\ A_j \end{pmatrix} = \begin{pmatrix} U_{11} & U_{12} \\ U_{21} & U_{22} \end{pmatrix} \begin{pmatrix} \widehat{R}_j \\ 0 \end{pmatrix} $
    \vspace{1mm}
    \State $A_{j+1} = U_{12}^HC_j$
    \State $B_{j+1} = U_{22}^HB_j$
    \State $l_{j+1} = l_j(a_jl_j^2 + b_j)/(c_jl_j^2 + 1)$ \Comment{Compute next value of $l$}
\EndFor
\State \Return $(A_p,B_p)$, optionally $l_p$
\end{algorithmic}
\end{algorithm}

\end{document}